\documentclass[12pt,a4paper]{amsart}
\usepackage{amssymb,enumerate}

\usepackage[pagebackref]{hyperref} 
\usepackage{longtable}

\newcommand{\FF}{{\mathbb{F}}}
\newcommand{\GG}{{\mathbb{G}}}
\newcommand{\LL}{{\mathbb{L}}}
\newcommand{\MM}{{\mathbb{M}}}
\newcommand{\QQ}{{\mathbb{Q}}}

\newcommand{\TT}{{\mathbb{T}}}
\newcommand{\ZZ}{{\mathbb{Z}}}

\newcommand{\bC} {\mathbf C}
\newcommand{\bG} {\mathbf G}
\newcommand{\bL} {\mathbf L}
\newcommand{\bT} {\mathbf T}
\newcommand{\bX} {\mathbf X}
\newcommand{\ba} {\mathbf a}
\newcommand{\bi} {\mathbf i}
\newcommand{\bu} {\mathbf u}

\newcommand{\cA} {\mathcal A}
\newcommand{\cE} {\mathcal E}
\newcommand{\cF} {\mathcal F}
\newcommand{\cG} {\mathcal G}
\newcommand{\cH} {\mathcal H}
\newcommand{\cM} {\mathcal M}
\newcommand{\cO} {\mathcal O}
\newcommand{\cU} {\mathcal U}
\newcommand{\cX} {\mathcal X}

\newcommand{\ff} {\mathfrak f}

\newcommand{\der}{{\operatorname{der}}}

\newcommand{\Sol}{{\operatorname{Sol}}}
\newcommand{\Frac}{{\operatorname{Frac}}}
\newcommand{\Ind}{{\operatorname{Ind}}}
\newcommand{\Irr}{{\operatorname{Irr}}}
\newcommand{\Res}{{\operatorname{Res}}}
\newcommand{\reg}{{\operatorname{reg}}}
\newcommand{\spets}{{\operatorname{spets}}}
\newcommand{\Spin}{{\operatorname{Spin}}}
\newcommand{\St}{{\operatorname{St}}}
\newcommand{\Uch}{{\operatorname{Uch}}}
\newcommand{\uni}{{\operatorname{uni}}}
\newcommand{\SL}{\operatorname{SL}}
\newcommand{\GL}{\operatorname{GL}}

\newcommand\RLG{{R_\bL^\bG}}

\newcommand{\tw}[1]{{}^{#1}\!}

\newcommand{\tu}{{\widetilde u}}
\newcommand{\tvhi}{{\widetilde\vhi}}
\newcommand{\wbu}{{\widetilde\bu}}
\newcommand{\bla}{{\big\langle}}
\newcommand{\bra}{{\big\rangle}}

\newcommand\Ph[1]{\Phi_{#1}}
\newcommand\hlf{\frac{1}{2}}
\newcommand\sq[1]{\sqrt{-#1}}
\newcommand\pl{{\!+\!}}
\newcommand\co{{:}}
\newcommand\mn{{\!-\!}}

\let\eps=\epsilon
\let\ga=\gamma
\let\La=\Lambda
\let\la=\lambda
\let\vhi=\varphi
\let\tht=\theta
\let\ze=\zeta

\newtheorem{thm}{Theorem}[section]
\newtheorem{lem}[thm]{Lemma}
\newtheorem{cor}[thm]{Corollary}
\newtheorem{prop}[thm]{Proposition}
\newtheorem{conj}[thm]{Conjecture}

\theoremstyle{definition}
\newtheorem{rem}[thm]{Remark}
\newtheorem{defn}[thm]{Definition}
\newtheorem{exmp}[thm]{Example}

\numberwithin{equation}{section}

\begin{document}

\title{Partial character tables for $\ZZ_\ell$-spetses}

\author{Radha Kessar}
\address{Department of Mathematics, University of Manchester, M13 9PL,
  United Kingdom}
\email{radha.kessar@manchester.ac.uk}

\author{Gunter Malle}
\address{FB Mathematik, RPTU Kaiserslautern, Postfach 3049,
  67653 Kaisers\-lautern, Germany.}
\email{malle@mathematik.uni-kl.de}

\author{Jason Semeraro}
\address{Department of Mathematical Sciences, Loughborough University, LE11 3TT,
  United Kingdom}
\email{j.p.semeraro@lboro.ac.uk}

\begin{abstract}
Let $\ell$ be a prime and $\GG$ be a simply connected $\ZZ_\ell$-spets, let
$q$ be a prime power, prime to $\ell$ and let $S$ be the underlying Sylow
$\ell$-subgroup of $\GG(q)$. Firstly,
motivated by known formulae for values of Deligne--Lusztig characters of
finite reductive groups, we propose a formula for the values of the
unipotent characters of $\GG(q)$ on the elements of $S$. Using this, we
explicitly compute the unipotent character values of the $\ZZ_2$-spets
$G_{24}(q)$ related to the Benson--Solomon fusion system $\Sol(q)$, and the
Aguad\'{e} $\ZZ_5$-spets $G_{29}(q)$.

Secondly, when $\ell>2$ is a very good prime for $\GG$, the Weyl group $W$ of
$\GG$ has order coprime with $\ell$ and $q\equiv1\pmod\ell$, we introduce a
formula for the values of characters in the principal block of $\GG(q)$ on
elements of $S$ which extends the Curtis--Schewe type formulae for groups of
Lie type, and which we show to satisfy a version of block orthogonality. 

In both cases we formulate and provide evidence for several conjectures
concerning the proposed values. 
\end{abstract}



\date{\today}

\maketitle

\pagestyle{myheadings}

\section{Introduction}
We continue our ongoing project to extend notions and results from the
$\ell$-modular representation theory of finite groups of Lie type in
cross-characteristic to spetses via the theory of $\ell$-compact groups.

Let $\ell$ be a prime and let $\cO\supset\ZZ_\ell$ be a complete discrete
valuation ring with field of fractions~$K$.
Recall from \cite[Def.~6.1]{KMS} that a pair $\GG=(W\phi,\La)$ consisting of a
$\ZZ_\ell$-lattice~$\La$ and a coset $W\phi$ with $W\le\GL(\La)$ a spetsial
$\ell$-adic reflection group and $\phi\in N_{\GL(\La)}(W)$ is called a
\emph{$\ZZ_\ell$-spets}.
Suppose in addition that $\GG$ is simply connected and $\ell$ is very good for
$\GG$ (\cite[Def.~2.4]{KMS}). Then for any prime power $q\in\ZZ_\ell^\times$,
in the papers \cite{KMS2,KMS} we attribute various data to the
pair $(\GG,q)$ including a Sylow $\ell$-subgroup $S$; a fusion system~$\cF$
on~$S$; a principal $\ell$-series $\cE_\ell(\GG(q),1)$ with degree function
$\deg:\cE_\ell(\GG(q),1)\to\QQ_\ell$; a set of characters in the
principal block $\Irr(B_0)\subseteq\cE_\ell(\GG(q),1)$ and an $\ell$-adic
integer $\dim(B_0)$. We furthermore show that these data recover the
corresponding notions for $\GG(q)$ in the `group case' where $W$ is rational
and $\GG(q)$ is the complete root datum of an actual finite group of Lie type.
We employ a variety of tools to demonstrate that the data satisfy familiar
group- or representation-theoretic properties. For example in \cite{KMS2},
the equality $\dim(B_0)_\ell=|S|_\ell$ (when $(|W|,\ell)=1$) is shown to follow
from the strong symmetry of the Hecke algebra $\cH(W,\bu)$ of $W$.

Suppose that $(\GG,q)$ is as above and let $\cF$ be the fusion system of
$\GG(q)$ over the $\ell$-group~$S$ (see \cite[Sect.~6.2]{KMS}). Here, we
initiate a programme of work whose goal is to extend the degree function $\deg$
to a \emph{value} function $\cX:\cE_\ell(\GG(q),1)\times S \to \cO$ which, when
it exists, we call the \emph{partial character table} of $\GG(q)$. In this case
we simply write $\chi(t):=\cX(\chi,t)$ to denote the value
$\chi\in\cE_\ell(\GG(q),1)$ takes at $t\in S$ so that $\cX(-,-)$ extends
$\deg(-)$ in the sense that $\deg(\chi)=\chi(1)$ for
$\chi\in\cE_\ell(\GG(q),1)$. These values should behave, as far as possible,
like those for actual finite groups. Three such desirable properties are:
\begin{itemize}
\item[(1)] viewed as an $S$-class function, $\chi\in\cE_\ell(\GG(q),1)$ is an
  $\cO$-linear combination of $\Irr(S)$;
\item[(2)] for $\chi\in\cE_\ell(\GG(q),1)$, the analogue of Frobenius' theorem
  holds (see Theorem \ref{t:frob});
\item[(3)] for $t,t'\in S$, the analogue of block orthogonality holds for
  $\Irr(B_0)$ (see Theorem \ref{t:centblock}). 
\end{itemize}

This paper advances the above programme in two directions. In
\S \ref{sec:spetschars}, under the additional hypotheses that $\phi=1$,
$(\ell,|W|)=1$ and $q\equiv1\pmod\ell$, we propose a function $\cX$ as above
with the property that, when $W$ is rational, $\chi(t)$ is the character value
of an actual split finite reductive group $G=\bG^F$ with Weyl group
$W=N_\bG(\bT)/\bT$ and maximally $F$-split torus $\bT$ (see
Proposition~\ref{p:blockbijection}) at the element $t\in\bT^F$. Our formula
for $\chi(t)$ is a natural generalisation of the Curtis--Schewe type formula
(see, for example, \cite[Cor.~3.3.13]{GM20}) for the character values of a
finite reductive group $\bG^F$ on its semisimple elements. For this
function, statements (1), (2) and (3) correspond respectively to
Conjectures~\ref{conj:restriction}, \ref{c:frobspets} and \ref{cor:main1}. We
emphasise that, when $W$ is rational, these conjectures follow from standard
group-theoretic results. We restate the last of the conjectures here:

\begin{conj}   \label{c:main1}
 Let $\GG=(W,\La)$ be a simply connected $\ZZ_\ell$-spets such that $\ell>2$ is
 very good for $\GG$ and coprime with $|W|$. Let $q\in1+\ell\ZZ_\ell$ and $B_0$
 be the principal block of $\GG(q)$. Then for all $t,t'\in S$,
 $$\sum_{\chi \in \Irr(B_0)} \chi(t)\chi(t'^{-1})
    = \begin{cases} \dim(B_0(t))& \text{if $t\sim_\cF t'$,}\\
                    0& \text{otherwise,}\end{cases}$$
 where $B_0(t)$ denotes the principal block of the $\ZZ_\ell$-spets $C_\GG(t)$.
\end{conj}

Here the $\ZZ_\ell$-spets $C_\GG(t)$ is defined in \cite[\S 6.2]{KMS} using the
theory of $\ell$-compact groups. We prove the following:

\begin{thm}   \label{t:main1}
 Conjecture~\ref{c:main1} holds for $\cX$ if $\cH(W,\bu)$ is strongly symmetric.
\end{thm}

According to \cite[Prop.~3.5]{KMS2}, all imprimitive irreducible spetsial
reflection groups apart possibly from $G(e,p,2)$ with $p$ even, as well as five
of the
primitive ones are known to have strongly symmetric Hecke algebras, and it is
conjectured that all of them do. Our strategy for proving Theorem~\ref{t:main1}
is similar to the approach taken in \cite{KMS2}. That is,
Conjecture~\ref{c:main1} follows from a more general statement,
Theorem~\ref{thm:B02}, for admissible 1-parameter specialisations of Hecke
algebras of $\ZZ_\ell$-reflection groups. 

In \S\ref{sec:valuesunip}, without the aforementioned additional hypotheses on
$\phi$, $|W|$ and $q$, and when $\ell$ is not necessarily very good, we give
formulae for the values of unipotent characters on $\ell$-elements. These are
shown to extend
those previously introduced when $\ell$ is very good, $(|W|,\ell)=1$, $\phi=1$,
and $q\equiv1\pmod\ell$. In this setting statements (1) and (2) above
correspond respectively to Conjectures~\ref{conj:restriction-uni}
and~\ref{c:frobspets-uni}. We moreover observe, using the associated Fourier
matrix, that each unipotent character $\chi$, associated to some
$\vhi\in\Irr(W\phi)$, lies in $\cO[\Irr(S)]$ if and only if for all
$\psi\in\Irr(S)$, the quantity
$$|S|^{-1}\sum_{t\in S}\psi(t)\, \big(f_\vhi^{C_\GG(t)}\big)_{x=q}$$
is a rational integer. Here $f_\vhi^\LL$ denotes the fake degree of
$\vhi|_{W_\LL}$ for $\LL$ (see, for example, \cite[\S1.2.2]{BMM14}). In
\S\ref{s:sol} we list the unipotent character values of the $\ZZ_\ell$-spets
$G_{24}(q)$ for the bad prime $\ell=2$. These are in turn used to supply some
evidence for
Conjectures~\ref{conj:restriction-uni} and~\ref{c:frobspets-uni}. A version of
statement~(3) would require formulae for values of non-unipotent character on
$\ell$-elements in this more general setting, and we defer this to future work. 
\medskip 

The paper is structured as follows. In \S\ref{sec:background} we record some
background results on character values and introduce various objects relating
to $\ZZ_\ell$-reflection cosets, extending those previously introduced in
\cite{KMS}.
In \S\ref{sec:chararbitrary} we define partial character tables for simply
connected $\ZZ_\ell$-reflection groups $W$ when $W$ is prime to $\ell$, and
prove that they are orthogonal provided that the underlying Hecke algebra is
strongly symmetric (Theorem~\ref{thm:B02}). In \S\ref{sec:spetschars} we define
the function $\cX(-,-)$ for
$q\in 1+\ell\ZZ_\ell$ (Definition \ref{def:chartab spets}) and apply
Theorem~\ref{thm:B02} to prove Theorem \ref{t:main1}, restated there as
Theorem~\ref{thm:blockorth}. Finally, in \S \ref{sec:valuesunip} we extend
Definition \ref{def:chartab spets} to give unipotent character values for
arbitrary spetses using Deligne--Lusztig characters. We give these values
explicitly for $G_{24}(q)$ and $G_{29}(q)$ respectively at the (bad) primes $2$ and $5$ when $q\in 1 + 4\ZZ_2$ and $q \in 1 + 5 \ZZ_5$
in~\S\ref{s:sol} and \S\ref{s:agu}.
\medskip

\noindent\textbf{Acknowledgements:}
The authors would like to thank Geoffrey R.~Robinson for providing illuminating
comments on and a reference for the proof of Theorem \ref{t:centblock}.
G.~Malle gratefully acknowledges financial support by the DFG--- Project-ID
286237555 -- TRR 195. J.~Semeraro gratefully acknowledges funding from the
UK Research Council EPSRC for the project EP/W028794/1. Finally, we sincerely
thank the anonymous referee for their careful and constructive reading of the
manuscript.

\section{Background}   \label{sec:background}

\subsection{On character values of finite groups} Here we introduce two general
character-theoretic statements which we will use later to motivate our
conjectures. Let $\ell$ be a prime number. 
The following result is a consequence of the second and third main theorems
of Brauer.

 \begin{thm}[Block orthogonality]   \label{t:centblock}
 Let $G$ be a finite group, $B$ be an $\ell$-block of $G$ and $s,s'\in G$ be
 $\ell$-elements. Let $C$ be the sum of the $\ell$-blocks of $C_G(x)$ which
 induce to $B$. Then
 $$\sum_{\chi\in\Irr(B)} \chi(s){\chi({s'}^{-1})}
   =\begin{cases} \dim(C) & \mbox{ if $s$ and $s'$ are $G$-conjugate; }\\
        0& \mbox{ otherwise. } \end{cases}$$
 Moreover, if $B$ is principal, then $C$ is the principal $\ell$-block of
 $C_G(x)$.
\end{thm} 

\begin{proof}
The first assertion is a consequence of Brauer's second main theorem and
follows from Corollary~8 and Theorem~5 of \cite{Rob}.  The second assertion is
Brauer's third main theorem (see for instance \cite[Thm~6.3.14]{LiBook2}).
\end{proof}

The next statement asserts that for any character $\chi$ of $G$, the partial
sum obtained from summands in the scalar product $\sum_{g\in G}\chi(g)$
which correspond to $\ell$-elements is divisible by $|G|_\ell$. 

\begin{thm}[Frobenius' Theorem]\label{t:frob}
 Let $G$ be a finite group, $\ell$ be a prime, $S$ be a Sylow $\ell$-subgroup
 of $G$ and $S^G$ be a set of $G$-class representatives of elements of $S$.
 For any character $\chi$ of $G$, we have
 $$\sum_{s \in S^G} |G:C_G(s)|\chi(s) \equiv 0 \pmod {|S|}.$$ 
\end{thm}

\begin{proof}
This follows from \cite[Cor.~41.9]{CR62} with $C$ the class of the trivial
element and $n=|G|_\ell=|S|$. Note that one may assume that $\chi$ is
irreducible.
\end{proof}

\subsection{$\ZZ_\ell$-reflection cosets}\label{subsec:zlref}
Let $\ell$ be a prime number. We rephrase some definitions on complex
reflection cosets from \cite[3A]{BMM99} to the $\ZZ_\ell$-setting. For this
recall that the \emph{parabolic subgroups} of a finite reflection group $W$ on
a vector space in characteristic~0 are defined to be the point-wise stabilisers
of subspaces. It is a result of Steinberg \cite[Thm~1.5]{St64} that these are
always reflection subgroups of $W$.

\begin{defn} Define the following:
\par
1. A \emph{$\ZZ_\ell$-reflection coset} $\GG$ is a pair $(W\phi,\La)$, where
  \begin{itemize}
    \item $\La$ is a $\ZZ_\ell$-lattice of finite rank,
    \item $W\le\GL(\La)$ is a finite reflection group, and
    \item $\phi\in\GL(\La)$ is an element of finite order normalising~$W$.
  \end{itemize}
 Note that only the coset $W\phi$ (and not the automorphism $\phi$) is defined
 by $\GG$. The rank of $\La$ is called the \emph{rank of $\GG$}.
 The $\ZZ_\ell$-reflection coset $\GG$ is \emph{split} if $W\phi = W$, i.e., if
 $\phi\in W$. If $\QQ_\ell\otimes\La$ is an irreducible $\QQ_\ell W$-module,
 then $\GG$ is called \emph{simple}.
\par
2. A \emph{sub-reflection coset} of $\GG = (W\phi,\La)$ is a
 $\ZZ_\ell$-reflection coset of the form $\GG' = (W'(w\phi)|_{\La'},\La')$,
 where $\La'$ is a pure sublattice of~$\La$, $W'$ is a reflection subgroup of
 $N_W(\La'){|_{\La'}}$ (the restriction to $\La'$ of the stabiliser $N_W(\La')$
 of $\La'$, isomorphic to $N_W(\La')/C_W(\La')$), and $w\phi\in W\phi$
 stabilises $\La'$ and normalises~$W'$.
\par
3. A \emph{toric $\ZZ_\ell$-reflection coset (or torus)} is a
 $\ZZ_\ell$-reflection coset with trivial reflection group. A \emph{torus of
 $\GG = (W\phi,\La)$} is a sub-coset of $\GG$ of the form
 $((w\phi)|_{\La'},\La')$ (so $\La'\le\La$ is a pure sublattice and
 $w\phi\in W\phi$ stabilises $\La'$).
\par
4. The \emph{centraliser} of a torus $\TT = ((w\phi)|_{\La'},\La')$ of $\GG$
 is the $\ZZ_\ell$-reflection coset $C_\GG(\TT):= (C_W(\La')w\phi,\La)$ (note
 that $C_W(\La')$ is a reflection subgroup of $W$ by Steinberg's theorem). Any
 such $\ZZ_\ell$-reflection coset is called a \emph{Levi} of~$\GG$. Thus, the
 Levis of $\GG$ are particular sub-reflection cosets of maximal rank of~$\GG$.
\par
5. The \emph{connected centre} of $\GG$ is the torus
 $Z^\circ(\GG):= (\phi|_{\La^W},\La^W)$ (notice that $\La^W$ is a pure
 sublattice of $\La$ and that the definition does not depend on the choice of
 $\phi$ in $W\phi$).
\end{defn}

Clearly any torus of $\GG$ is contained in some \emph{maximal} torus of $\GG$, 
that is, a $\ZZ_\ell$-reflection coset of the form $(w\phi,\La)$ where $w\in W$.

\begin{lem}   \label{lem:levis}
The following hold:
 \begin{enumerate}[\rm(a)]
  \item The Levis of $\GG = (W\phi,\La)$ are the $\ZZ_\ell$-reflection cosets
   $\LL= (W'w\phi,\La)$, where $W'\le W$ is a parabolic subgroup and
   $w\phi\in W\phi$ normalises $W'$.
  \item If $\LL$ is a Levi of $\GG$, then $C_\GG(Z^\circ(\LL)) = \LL$.
 \end{enumerate}
\end{lem}

\begin{proof}
(See \cite[Lemma 3.2]{BMM99}.) By the above definition of Levis it suffices to
prove (b) for any $\LL$ as in (a). Set $\La':=\La^{W'}$. By definition
$Z^\circ(\LL) = ((w\phi)_{|_{\La'}},\La')$; hence
$C_\GG(Z^\circ(\LL)) = (C_W(\La')w\phi,\La)$, and the result follows from
Steinberg's theorem.
\end{proof}

Suppose that $(W,\La)$ is a $\ZZ_\ell$-reflection group. A reflection subgroup
$W_0\le W$ is \emph{full}, if whenever $r\in W$ is a reflection and
$1\ne r^n\in W_0$ then $r\in W_0$. In particular all parabolic subgroups
of~$W$ are full reflection subgroups of~$W$. The following observation will be
needed later.

\begin{lem}   \label{lem:stab}
 If $(W,\La)$ is a $\ZZ_\ell$-reflection group with $\ell\nmid|W|$, then for
 any $a\ge 1$ and $s\in \La/\ell^a \La$, $C_W(s)$ is parabolic and hence a
 full reflection subgroup of~$W$.
\end{lem}

\begin{proof}
This follows from \cite[Prop.~2.3]{KMS}.
\end{proof}

\begin{defn}   \label{def:2.5}
 Recall that a $\ZZ_\ell$-reflection coset $\GG=(W\phi,\La)$ is said to be
 \emph{simply connected} if $\La=[\La,W]$, so if its \emph{fundamental group}
 $\La/[\La,W]$ is trivial.
 We remark that when $\ell$ is odd, this is equivalent to $\Omega X$ being
 simply connected where $X$ is the $\ell$-compact group associated to $(W,\La)$
 (see \cite[Sec.~3.2]{KMS}).

 For $w\in W$ and $H\le W$ a $w\phi$-stable subgroup we denote
 $(Hw\phi,\La)^\circ:=(H^\circ w\phi,\La)$, where $H^\circ$ is the subgroup
 of~$H$ generated by its reflections, obviously $w\phi$-stable as well.
\end{defn}

Following \cite[Def.~1.44(nc)]{BMM14} the \emph{order (polynomial)} of a
$\ZZ_\ell$-reflection coset $\GG=(W\phi,\La)$ is defined as
$$|\GG|=|\GG|(x):=x^{N(W)}\prod_{i=1}^r \eps_i^{-2}(x^{d_i}-\eps_i)
  \in \ZZ_\ell[x]$$
where $N(W)$ is the number of reflections in $W$ and
$\{(d_i,\eps_i)\mid 1\le i\le r \}$ is the uniquely determined multiset
of generalised degrees of $W\phi$ for its action on $\La\otimes\QQ_\ell$. (This
differs by a root of unity from the order polynomial defined in
\cite[(3.4)]{BMM99}.) Let $\Phi\in\ZZ_\ell[x]$ be an irreducible cyclotomic
polynomial. A torus $\TT$ is a \emph{$\Phi$-torus} if $|\TT|$ is a power
of~$\Phi$. The centralisers of $\Phi$-tori of $\GG$ are called
\emph{$\Phi$-split} Levis. Since $\Phi$-tori are defined in terms of
eigenspaces in $\La\otimes\overline{\QQ_\ell}$, the Sylow theorems
for maximal $\Phi$-tori over $\QQ_\ell$ (see \cite[Thm~1.50]{BMM14}) continue
to hold over $\ZZ_\ell$; these maximal $\Phi$-tori are hence called \emph{Sylow
$\Phi$-tori} of $\GG$.
For $\LL=(W_\LL w\phi,\La)$ a $\Phi$-split Levi of $\GG$ we define its
\emph{relative Weyl group} $W_\GG(\LL):=N_W(W_\LL w\phi)/W_\LL$. 

Now assume that $W$ is a \emph{spetsial} $\ZZ_\ell$-reflection group on a
$\ZZ_\ell$-lattice $\La$. See \cite[\S3]{MaICM} for the definition and
classification of complex spetsial reflection groups; the $\ZZ_\ell$-reflection
groups (among them) can be deduced for example from \cite[\S2.1]{KMS}. In
this case any reflection coset $\GG=(W\phi,\La)$ is called a
\emph{$\ZZ_\ell$-spets} (see e.g.\ \cite[Def.~6.1]{KMS}).

\subsection{Fusion systems and centralisers}   \label{subsec:centr}
Let $\GG=(W\phi,\La)$ be a $\ZZ_\ell$-reflection coset where
$\phi\in N_{\GL(\La)}(W)$ is of $\ell'$-order, and assume $(W,\La)$ is simply
connected 

Let $\ze\in\ZZ_\ell^\times$ be a root of unity and $a\ge1$. If $\ell $ is odd,
then via the theory of $\ell$-compact groups, attached to these data plus a
choice of $q\in\ze+\ell^a\ZZ_\ell^\times$ there comes a fusion system~$\cF$ on
a finite $\ell$-group~$S$ (see~\cite[Thm~3.2]{KMS}). We say $S$ is the
\emph{Sylow $\ell$-subgroup of~$(\GG,q)$}. Here, $S$ is an extension of a
homocyclic $\ell$-group $T$ of exponent $\ell^a$ (the \emph{toral part}) by a
Sylow $\ell$-subgroup of the associated relative Weyl group $W_{\phi\ze^{-1}}$
(see \cite[Thm~3.6]{KMS}). For $s\in S$ let $W(s)$ and $\phi_s$ be as defined
after \cite[Prop.~5.5]{KMS}; recall that $\phi_s\in W\phi$ normalises $W(s)$.
Thus, for any $s\in S$ this defines a $\ZZ_\ell$-sub-reflection coset
$C_\GG(s):=(W(s)\phi_s,\La)$ of~$\GG$. Note $C_{\GG}(s)$ is well defined (see
the remarks preceding \cite[Lemma 6.6]{KMS}).

\begin{lem}   \label{lem:phi_t}
 Suppose that $\ell>2$, $\ze=1$, $\phi=1$ and $\ell\nmid|W|$. Then $S=T$ may
 be identified with $\La/\ell^a \La$ and under the inherited action of $W$ on
 $\La/\ell^a \La$, $\cF=\cF_\ell(T W)$.
 Moreover, for any $s\in S$ we have $W(s) = C_W(s)$ and $\phi_s=1$.
\end{lem}

\begin{proof}
The first three assertions are shown in \cite[Ex.~4.4]{KMS2}. Note that the
spetsial assumption in \cite[Ex.~4.4]{KMS2} is not relevant for the arguments.
The other assertions are implicit in \cite[Ex.~4.4]{KMS2} and follow from the
explicit descriptions of $W(s)$ and $\phi_s$ after the proof of
\cite[Prop.~5.5]{KMS}.
\end{proof}

\section{Deforming character tables of Sylow normalisers}   \label{sec:chararbitrary}
Throughout this section $W$ is an arbitrary (finite) $\ZZ_\ell$-reflection
group on a $\ZZ_\ell$-lattice~$\La$. We assume moreover that $|W|$ is prime
to~$\ell$. Let $\bu=(u_{r,i})$ be a set of indeterminates where $r$ runs over
conjugacy class representatives of distinguished reflections of $W$ and
$0\le i<o(r)$, and let $\cH(W,\bu)$ be the \emph{(generic) Hecke algebra} over
$\ZZ[\bu^{\pm1}]$ associated
to $W$ (see \cite[Def.~4.21]{BMR}). For certain specialisations of $\cH(W,\bu)$
we introduce and study deformations of the character table of the associated
torus normaliser that will subsequently be used to define partial character
tables of $\ZZ_\ell$-spetses.

\subsection{The principal block}   \label{subsec:B0 refl}
Let $\wbu=(\tu_{rj})$ be indeterminates such that
$$\tu_{rj}^z=\exp(-2\pi\bi j/o(r))u_{rj}\qquad\text{for all $r,j$}$$
where $z:=|Z(W)|$. Thus, $K:=\Frac(\ZZ_\ell[\wbu])$ is a splitting field
for $\cH(W,\bu)$ by \cite[Thm~5.2]{Ma99}. Let $R$ be an integral domain
containing $\ZZ_\ell$ and $x^{1/z}$ an indeterminate. An \emph{admissible
specialisation} of $\cH(W, \bu)$ is one induced by an $R$-linear ring
homomorphism of the form
$$\psi_\ba:R[\wbu^{\pm1}]\to R[x^{\pm1/z}],\quad
  \widetilde u_{rj}\mapsto x^{a_{rj}/z}\qquad\text{for all $r,j$},
\eqno{(\ddagger)}$$
for a tuple of integers $\ba=(a_{rj})$ (here specialisation is in the sense of
\cite[\S3B]{KMS2}, applied with $R[\wbu^{\pm1/z}]$ in place of $A'$ and
$R[x^{\pm1/z}]$ in place of $R$). Observe that if
$\psi_\ba$ is admissible, its composition with the specialisation $x\mapsto1$
is the canonical specialisation to the group algebra. Also, recall
the Schur elements $\ff_\vhi\in \ZZ_\ell[\wbu^{\pm 1}]$ of $\cH(W,\bu)$ for
$\vhi\in\Irr(W)$ (see e.g.~\cite[\S3A]{KMS2}). 


For any integer $a\ge1$ and admissible specialisation $\psi_\ba$ for
$\cH(W,\bu)$ we write $\GG_{\ba,a}:=(W,\La,\psi_\ba,a)$.  Set 
$T=\La/\ell^a \La$ and for $s\in T$, define the centraliser of $s\in T$ to
be $C_{\GG_{\ba,a}}(s):=(C_W(s),\La)$. Note that by Lemma~\ref{lem:stab},
$W_0:=C_W(s)$ is a full reflection subgroup of $W$. 
We denote by $\cH(W_0,\bu_0)$ the Hecke algebra of
$W_0$ whose parameters $\bu_0$ consist of those parameters for~$W$ whose
corresponding distinguished reflections are, up to conjugation, contained
in~$W_0$. It follows from the Freeness Theorem (see \cite[Thm~3.1]{KMS2} and
the references given there) that $\cH(W_0,\bu_0)$ is naturally a subalgebra of
$\cH(W,\bu)$, for which $K$ is still a splitting field (see the
explicit results in \cite[Thm~5.2]{Ma99}). Thus, an admissible specialisation
$\psi_\ba$ for $\cH(W,\bu)$ gives rise to an admissible specialisation
$\psi_{s,\ba}$ for $\cH(C_W(s),\bu_0)$ as follows: we define $\psi_{s,\ba}$ as
the restriction of $\psi_\ba$ to the parameters associated to the distinguished
reflections in $C_W(s)$. It is clear that $\psi_{s,\ba}$ is admissible.
 
Since $|W|$ is prime to $\ell$, Gallagher's theorem gives
$$\Irr(TW)=\{\rho_{\tht,\vhi}\mid\tht\in\Irr(T)/W,\ \vhi\in\Irr(W(\tht))\}$$
where $W(\tht):=C_W(\tht)$ and
$$\rho_{\tht,\vhi} = \Ind_{TW(\tht)}^{TW}(\tht\otimes\vhi)$$
(by abuse of notation denoting by $\tht$ the trivial extension of the linear
character $\tht$ of $T$ to $TW(\tht)$ and by $\vhi$ the inflation of $\vhi$
to the same group). This is the blueprint for our definition of the principal
block of $\GG_{\ba,a}$, generalising the definition of \cite[Def.~4.1]{KMS2}
(see Section~\ref{subsec:chartab spets}).
For this we choose and fix a $W$-equivariant bijection
$\widehat{\ }: T\to\Irr(T)$ (see \cite[\S4B]{KMS2}).

\begin{defn}   \label{def:B}
 Recall that $|W|$ is prime to $\ell$. Let $\psi_\ba$ be an admissible
 specialisation for $\cH(W,\bu)$ and $a\ge1$. To each $\tht\in\Irr(T)/W$ and
 $\mu\in\Irr(W(\tht))$ we associate a symbol $\ga_{\tht,\mu}$ and refer
 to the set
 $$\Irr(B_0):=
   \{\ga_{\tht,\mu}\mid\tht\in\Irr(T)/W,\ \mu\in\Irr(W(\tht))\}$$
 as \textit{the characters in the principal block $B_0=B_0(\GG_{\ba,a})$ of
 $\GG_{\ba,a}$}. This set comes equipped with a map $\deg:\Irr(B_0)\to 
 \QQ_\ell(x^{1/z})$ defined as
 $$\deg(\ga_{\tht,\mu}):=\psi_\ba(\ff_{1_W})/\psi_{s,\ba}(\ff_\mu)
   \ \in\ \QQ_\ell(x^{1/z})$$
 where $s \in T$ is such that $\tht=\widehat{s}$. Note that we also write
 $\ga_{\tht,\mu}(1)$ for $\deg(\ga_{\tht,\mu})$ and refer to this element as
 the \textit{degree} of $\ga_{\tht,\mu}$.

 We moreover set $\dim(B_0):=\sum_{\ga\in\Irr(B_0)}\ga(1)^2$. Note that since
 the canonical specialisation of $\cH(C_W(s),\bu_0)$ to the group algebra of
 $C_W(s)$ factors through $\psi_{s,\ba}$, and since the image of the Schur
 element $\ff_\mu$ under the canonical specialisation is the reciprocal of the
 codegree of $\mu$, $\psi_{s,\ba}(\ff_\mu)$ is non-zero.
 \end{defn}

\subsection{Partial character tables for $\ZZ_\ell$-reflection groups}
As we will be dealing with characters in principal blocks of different
reflection data, we write $\ga_{\tht,\vhi}^W$ for $\ga_{\tht,\vhi}$ if we want
to specify that the character is with respect to $W$. Now fix an extension
field~$F$ of $\QQ_\ell$ containing a primitive $\ell^a$-th root of unity and
identify the characters of $T$ with $F$-valued characters of~$T$.

\begin{defn}   \label{def:chartab}
 For $\ga_{\tht,\vhi}^W\in\Irr(B_0)$ (with $\tht\in\Irr(T)$ and
 $\vhi\in\Irr(W(\tht))$) and $t\in T$ with $W_t:=C_W(t)$, set
 $$\ga_{\tht,\vhi}^W(t):=
   \sum_{w\in W_t\backslash W/W(\tht)}\tht(t^w)\sum_{\la\in \Irr(W_{t^w}(\tht))}
   \bla\vhi,\Ind_{W_{t^w}(\tht)}^{W(\tht)}(\la)\bra\,\ga_{\tht,\la}^{W_{t^w}}(1) \in F(x^{1/z})
 $$
 where $W(\tht):=C_W(\tht)$ and $W_t(\tht):=C_{W_t}(\tht)=C_W(\tht,t)$ 
 (note that by Lemma~\ref{lem:stab} these are parabolic subgroups of $W$).
 The matrix $\cX=\cX(B_0)=(\ga_{\tht,\vhi}^W(t))$ with rows indexed by $\Irr(B_0)$ and columns indexed by $T/W$ is
 called the \emph{(partial) character table of $B_0(\GG_{\ba,a})$}.
\end{defn} 

Note that, since
$\ga_{\tht,\la^{\nu}}^{W_{t^{\nu}}}(1)=\ga_{\tht,\la}^{W_t}(1)$ for any
$t\in T$, $\tht\in\Irr(T),\nu\in W(\tht)$ and $\la\in\Irr(W_t)$, the
definition of $\ga_{\tht,\vhi}^W(t)$ depends only on the $W$-class of
$(\tht,\vhi)$. In particular, $\ga_{\tht,\vhi}^W$ is $W$-invariant.

\begin{exmp}   \label{ex:central}
When $W_t=W$, that is, $t$ is centralised by $W$, this yields
$\ga_{\tht,\vhi}^W(t)=\ga_{\tht,\vhi}^W(1)\,{\tht}(t)$. At the
other extreme, if $t$ is regular, so $C_W(t)=1$, then
$$\ga_{\tht,\vhi}^W(t)=\vhi(1)\sum_{w\in W/W(\tht)}\tht(t^w)$$
as $\ga_{\tht,1}^1(1)=1$ for all $\tht$, by definition.
\end{exmp}

As announced above, the partial character table is a deformation of the columns
corresponding to $\ell$-elements of the character table of the Sylow
normaliser $TW$:

\begin{prop}   \label{prop:spec}
 Under the specialisation $x\mapsto1$, the character table $\cX(B_0)$ maps to
 the character table of $TW$, restricted to the columns for elements of $T$.
\end{prop}

\begin{proof}
Recall the identification
$$\Irr(TW)=\{\rho_{\tht,\vhi}\mid \tht\in \Irr(T)/W,\ \vhi\in\Irr(W(\tht))\},$$
where $W(\tht):=C_W(\tht)$, stated above. Then the natural map
$$\Theta:\Irr(B_0)\to\Irr(TW),\qquad \mbox{ defined by } \qquad
  \ga_{\tht,\vhi}^W \mapsto\rho_{\tht,\vhi},$$
is a bijection with $\Theta(\chi)(1)=\chi(1)|_{x=1}$. Indeed, since $\psi_\ba$
is admissible, its composition with the specialisation $x\mapsto1$ is the
canonical specialisation to the group algebra, and the same is true for
$\psi_{s,\ba}$ for $s\in T$. So, for any $\tht\in\Irr(T)$,
$\vhi\in\Irr(W(\tht))$, and $s\in T$ such that $\tht=\widehat{s}$ we have 
$$(\psi_\ba(\ff_{1_W})/\psi_{s,\ba}(\ff_\mu))_{x=1}= (|W| \vhi(1))/ |C_W(s)|
  = \rho_{\tht,\vhi}(1). $$
Here the first equality follows from \cite[Lemma~3.6]{KMS2} and the second is
just the formula for the degree of an induced character.
Now for $t\in T$ we have
$$\begin{aligned}
 \Theta(\ga_{\tht,\vhi}^W)(t)=& \rho_{\tht,\vhi}(t)
  = \frac{1}{|TW(\tht)|}\sum_{g\in TW}\vhi(1)\,{\tht}(t^g)\\
  =& \frac{1}{|W(\tht)|}\sum_{w\in W}\vhi(1)\,{\tht}(t^w)
  = \sum_{w\in W_t\backslash W/W(\tht)}
    \vhi(1)|W_{t^w}:W_{t^w}(\tht)|\,{\tht}(t^w).
\end{aligned}$$
Using that
$$\vhi|_{W_{t^w}(\tht)}=\sum_{\la\in\Irr(W_{t^w}(\tht))}\bla\vhi|_{W_{t^w}(\tht)},\la\bra\,\la
  =\sum_{\la\in\Irr(W_{t^w}(\tht))}\bla\vhi,\Ind_{W_{t^w}(\tht)}^{W(\tht)}(\la)\bra\,\la,$$
we find
$$\Theta(\ga_{\tht,\vhi}^W)(t)
  =\sum_{w\in W_t\backslash W/W(\tht)}\!\!\!\tw{w}{\tht}(t)|W_{t^w}:W_{t^w}(\tht)|
   \sum_{\la\in\Irr(W_{t^w}(\tht))}
   \!\!\!\bla\vhi,\Ind_{W_{t^w}(\tht)}^{W(\tht)}(\la)\bra\,\la(1)
  = \ga_{\tht,\vhi}^W(t)|_{x=1}$$
as claimed, since
$\ga_{\tht,\la}^{W_{t^w}}(1)|_{x=1}=|W_{t^w}:W_{t^w}(\tht)|\,\la(1)$
(see above).
\end{proof}

\subsection{Block orthogonality}
Motivated by Theorem \ref{t:centblock}, we now show that characters in the
principal block $B_0$ of $\GG_{\ba,a}$ are orthogonal in the following sense,
at least when the Hecke algebra $\cH(W,\bu)$ is strongly symmetric (see
also Remark~\ref{rem:part(2)} below).

\begin{thm}   \label{thm:B02}
 Let $\GG=(W,\La)$ be a simply connected $\ZZ_\ell$-reflection group with $|W|$
 prime to~$\ell$, $\psi_\ba$ an admissible specialisation, $a\ge1$, $T$ the
 Sylow $\ell$-subgroup of $\GG_{\ba,a}$ and $B_0$ the principal block
 of~$\GG_{\ba,a}$. If $\cH(W,\bu)$ is strongly symmetric, then for any
 $t,t' \in T$,
 $$\sum_{\chi \in \Irr(B_0)} \chi(t){\chi({t'}^{-1})}
   =\begin{cases} \dim(B_0(t)) & \mbox{ if $t$ and $t'$ are $W$-conjugate}\\
                             0 & \mbox{ otherwise,} \end{cases}$$
 where $B_0(t)$ denotes the principal block of $C_\GG(t)_{\ba,a}$ (with respect
 to the admissible specialisation $\psi_{t,\ba}$ for $\cH(C_W(t),\bu_0)$
 as above).
\end{thm}

We give the proof in steps. Note that the claim is tautological in
the case $t=t'=1$.

\begin{lem}   \label{lem:fixparabolic}
 Suppose that $W$ is an $\ell'$-group. Let $1\ne t\in T$ and let $H \leq W_t$.
 Let $A_H\leq \widehat T$ be the subset consisting of all characters $\tht$
 such that $W_t(\tht)$ is $W$-conjugate to $H$. Then $A_H$ is $W_t$-invariant
 and $\sum_{\tht\in A_H/W_t}\tht(t) =0$.
\end{lem}

\begin{proof}
The first assertion is immediate. Let $X=\{\tht\in\widehat T\mid\tht(t) =1\}$. 
Then $X$ is a $W_t$-invariant subgroup of $\widehat T$. Since $W_t$ is an
$\ell'$-group and $\widehat T$ is an $\ell$-group there is a $W_t$-invariant
decomposition $\widehat T = X' \times X$. On the other hand, if $\tht\in X'$,
then $\tht X$ is also $W_t$-invariant since for any $w\in W_t$,
$\tw{w}\tht(t) = \tht(\tw{w}t)=\tht (t)$. It follows that $W_t$ acts
trivially on $X'$. Hence, $W_t(\rho) =W_t(\rho\tht)$ for all
$\rho\in X$, $\tht\in X'$. It follows that $A_H = X' \times (A_H \cap X)$
and $\rho\mapsto\tht\rho$ induces a bijection between the set of $W_t$-orbits
on $X$ and the set of $W_t$-orbits on $\tht X$. Since
$\tht\rho (t) =\tht(t)$ for all $\rho\in X$, $\tht\in X'$ and since
$\sum_{\tht\in X'}\tht(t) =0$, we obtain
$$\sum_{\tht\in A_H/W_t}\tht(t)
   = |(A_H \cap X)/W_t|\sum_{\tht\in X'} \tht(t) = 0.\qedhere $$
\end{proof}

Let $\tau = t_{W,\bu} = \sum_{\vhi\in\Irr(W)}(f^W_{\vhi})^{-1}\vhi$ be the
trace form on $KW$ inherited from the canonical trace form on
$\cH(W,\bu)$ (cf. \cite[Conj.~3.3]{KMS2}) and set
$g^W_\vhi:= \psi_\ba(f^W_\vhi)^{-1}$. If $\cH(W,\bu)$
is strongly symmetric (see \cite[Def.~3.4]{KMS2}), then for any parabolic
subgroup $W_0$ of $W$, $\tau|_{KW_0}$ is the trace form
$\sum_{\la\in\Irr(W_0)}(\ff_\la^{W_0})^{-1}\la$ on $KW_0$ and by
abuse of notation we will denote this restriction also by $\tau$, as well its
specialisation under $\psi_\ba$.

For parabolic subgroups $W_0 \leq W_1$ of $W$ and for $\eta\in\Irr(W_1)$ denote
by $\ff_{W_0,\eta}^{W_1}$ the Schur element of the induced form
$\Ind_{W_0}^{W_1}\tau$ at $\eta$ and set
$g_{W_0,\eta}^{W_1}:= \psi_\ba(\ff_{W_0,\eta}^{W_1})^{-1}$.
For $\tht \in \Irr(T)$ and $t,t'\in T$ define
$$\Phi_\tht(t,t'):=
 \sum_{\vhi\in\Irr(W(\tht))}\ga_{\tht,\vhi}^W(t)\ga_{\tht,\vhi}^W({t'}^{-1})
  \qquad\text{and}\qquad a_{t,t'}:=\psi_\ba(\ff_{1_{W_t}}\ff_{1_{W_{t'}}}).$$

\begin{prop}   \label{prop:induced}
 Assume that $\cH(W,\bu)$ is strongly symmetric. Then for all $\tht\in\Irr(T)$
 and $t,t'\in T$,
 $$\Phi_\tht(t,t') = a_{t,t'}\sum_{w, w'}\tht(t^w ({t'}^{w'})^{-1}) 
   \ \sum_y \sum_{\la\in\Irr(W_{t^w}(\tht))} g_\la^{W_{t^w}(\tht)} g_{W_{t^w}(\tht)\cap W_{\!{t'}^{w'y^{-1}}}(\tht), \la}^{W_{t^w}(\tht)}$$
 where $w$ ranges over $W_t\backslash W/ W(\tht)$, $w'$ over
 $W_{t'}\backslash W/ W(\tht)$ and $y$ over
 $W_{t^w}(\tht)\backslash W(\tht)/W_{{t'}^{w'}}(\tht)$.
\end{prop}

\begin{proof}
We have
$$\Phi_\tht(t,t') =
  a_{t,t'}\sum_{w, w'} \sum_{\la\in\Irr(W_{t^w}(\tht))} g_\la^{W_{t^w}(\tht)} \!\!\sum_{\vhi\in\Irr(W(\tht))} \bla\vhi,\Ind_{W_{t^w}(\tht)}^{W(\tht)}(\la)\bra \,h(w',\vhi)$$
with $h(w',\vhi):=\sum_{\la'\in\Irr(W_{{t'}^{w'}}(\tht))} g_{\la'}^{W_{{t'}^{w'}}(\tht)} \bla\vhi,\Ind_{W_{{t'}^{w'}}(\tht)}^{W(\tht)}(\la')\bra$.
Now since $\cH(W,\bu)$ is strongly symmetric, $h(w',\vhi)$ is the coefficient
of $\vhi$ in $\Ind_{W_{{t'}^{w'}}(\tht)}^{W(s)} \tau$ (see the commentary
above the statement of this proposition), hence
$$\sum_{\vhi\in\Irr(W(\tht))}\bla\vhi,\Ind_{W_{t^w}(\tht)}^{W(\tht)}(\la)\bra\,
  h(w',\vhi) $$
is the coefficient of $\la $ in $\Res^{W(\tht)}_{W_{t^w}(\tht)} \Ind_{W_{{t'}^{w'}}(\tht)}^{W(\tht)} \tau$. By the Mackey formula,
$$\Res^{W(\tht)}_{W_{t^w}(\tht)} \Ind_{W_{{t'}^{w'}}(\tht)}^{W(\tht)} \tau = \sum_{y\in W_{t^w}(\tht)\backslash W(\tht)/W_{{t'}^{w'}}(\tht)} \Ind_{W_{t^w}(\tht)\cap\tw{y} W_{{t'}^{w'}}(\tht)} ^{W_{t^w}} \Res^{ \tw{y} W_{{t'}^{w'}}(\tht) }_{ W_{t^w}(\tht) \cap \tw{y}W_{{t'}^{w'}}(\tht)} \tw{y} \tau.$$ 
The result follows since by the strongly symmetric property for all
$y\in W(\tht)$,
$$\Res^{ \tw{y} W_{{t'}^{w'}}(\tht) }_{ W_{t^w}(\tht) \cap \tw{y}W_{{t'}^{w'}}(\tht)} \tw{y} \tau = \tau_{ W_{t^w}(\tht) \cap \tw{y}W_{{t'}^{w'}}(\tht)}\quad\text{and}\quad
\tw{y}W_{{t'}^{w'}}(\tht) = W_{{t'}^{w'y^{-1}}}(\tht).\qedhere$$
\end{proof} 

\begin{prop}   \label{prop:t1}
 Theorem~$\ref{thm:B02}$ holds if either $t$ or $t'$ is central.
\end{prop} 

\begin{proof}
Without loss of generality we may assume $t'$ is central. For a parabolic
subgroup $W_0$ of $W$ set
$\Phi^{W_0}(t,t'):= \sum_\chi\chi (t)\chi({t'}^{-1})$, where $\chi$ ranges over
the irreducible characters of the principal block of $(W_0,\La)$ and set
$\Phi(t,t'):= \Phi^{W}(t,t')$.
Note that $\Phi(t,t')= \sum_{\tht\in\widehat{T}/W} \Phi_\tht(t,t')$. It follows
from the definition of $\chi(t)$ that $\Phi(t,t') =\Phi(t{t'}^{-1},1)$, hence
we may also assume that $t'=1$ and $t\ne 1$. By Proposition~\ref{prop:induced},
$$\Phi(t,1)
 = a_{t,1}\sum_{\tht\in \widehat{T}/W}\sum_{w\in W_t\backslash W/W(\tht)}
  \tht(t^w)\sum_{\la\in \Irr(W_{t^w}(\tht))} (\ga_{\tht,\la}^{W_{t^w}}(1))^2.$$
Let $I$ be a set of representatives of $\widehat T/W $ and
for each $\tht \in I $, let $J_{\tht} $ be a set of representatives of
$W_t\backslash W/W(\tht) $. 
For a subgroup $H$ of $W_t$, let
$$\cA_H =\{ (\tht, w)\mid\tht\in I, w\in J_\tht,\,
  W_t(\tw{w}\tht) \text{ is }W\text{-conjugate to } H\}.$$
If $\tht_1,\tht_2\in \widehat{T}$ and $ w_1,w_2\in W$ are such that
$W_{t^{w_1}}(\tht_1)$ and $W_{t^{w_2}}(\tht_2)$ are $W$-conjugate, then 
$$ \sum_{\la\in \Irr(W_{t^{w_1}}(\tht_1))} (\ga_{\tht_1,\la}^{W_{t^{w_1}}}(1))^2
  = \sum_{\la\in \Irr(W_{t^{w_2}}(\tht_2))}(\ga_{\tht_2,\la}^{W_{t^{w_2}}}(1))^2.$$ 
Hence in order to prove the proposition (and noting the duality between $T$ and
$\widehat T$), it suffices to prove that 
$$\sum_{(\tht, w) \in \cA_H} \tw{w}\tht(t) =0 $$
for all $ H \leq W_t $.

Let $A_H$ be as in Lemma~\ref{lem:fixparabolic}. The map
$(\tht, w)\mapsto \tw{w}\tht$ induces a bijection between $\cA_H$ and
$A_H/W_t$ and now the result follows from the above displayed equation
for~$\Phi(t, 1)$.
\end{proof}

\begin{proof}[Proof of Theorem~\ref{thm:B02}]
Let $\tht\in \widehat{T}$. For $w,w'\in W$ we have
$$\tht( t^w ({t'}^{w'})^{-1}) = \eta (t ({t'}^{w' w^{-1}})^{-1})$$
where $\eta =\tw{w}{\tht} $. Now, $W_t(\eta) = \tw{w}W_{t^w}(\tht)$, $W_{\,{t'}^{(w'w^{-1})\tw{w}y^{-1}}} (\eta)= \tw{w}W_{\!{t'}^{w'y^{-1}}}(\tht)$ 
and for $\la\in\Irr(W_{t^w}(\tht))$, we have 
$$g_{\tw{w}\!\la}^{\tw{w}(W_{t^w}(\tht))} =g_\la^{W_{t^w}(\tht)}\quad\text{and}\quad
  g_{\tw{w}(W_{t^w}(\tht)\cap W_{\!{t'}^{w'y^{-1}}}(\tht)),\tw{w} \la}^{\tw{w}(W_{t^w}(\tht))}= g_{W_{t^w}(\tht)\cap W_{\!{t'}^{w'y^{-1}}}(\tht), \la}^{W_{t^w} (\tht)}.$$
Hence, replacing $\tht$ by $\eta$, $w'$ by $w'w^{-1}$, $y$ by $\tw{w}y$ and
$\la$ by $\tw{w}\la$ in Proposition~\ref{prop:induced}, we obtain
$$\Phi_\tht(t,t') = a_{t,t'}\sum_{\eta,w'}\sum_y\ \eta(t({t'}^{w'})^{-1}) 
 \!\sum_{\la\in\Irr(W_t(\eta))} g_\la^{W_{t}(\eta)} g_{W_t(\eta)\cap W_{\!{t'}^{w'y^{-1}}}(\eta),\la}^{W_t(\eta)},$$
where $\eta$ ranges over a set of representatives of the $W_t$-orbits on the
$W$-orbit of $\tht$, for each $\eta$, $w'$ ranges over a set of
representatives of $W_{t'} \backslash W/W(\eta)$ and for each $w'$,
$y$ ranges over a set of representatives of
$W_t(\eta)\backslash W(\eta)/W_{{t'}^{w'}}(\eta)$.
 
For a fixed $\eta$, as $w'$ ranges over a set of representatives of
$W_{t'}\backslash W /W(\eta)$ and $y$ ranges over a set of representatives of
$W_t (\eta)\backslash W(\eta)/W_{{t'}^{w'}}(\eta)$, $w'y^{-1}$ runs over
a set of representatives of $W_{t'}\backslash W/ W_t(\eta)$. Since
$\eta(u^y ) = \eta(u)$ for $y\in W(\eta)$ and $u\in T$, we obtain 
$$\Phi_\tht(t,t') = a_{t,t'}\sum_\eta\sum_{w'}\ \eta(t ({t'}^{w'})^{-1}) 
  \,h(\eta,w'),$$
with
$$h(\eta,w'):=\sum_{\la\in\Irr(W_t(\eta))} g_\la^{W_t(\eta)} g_{W_t(\eta)\cap W_{\!{t'}^{w'}}(\eta),\la}^{W_t(\eta)},$$
where $\eta$ ranges over a set of representatives of the $W_t$-orbits on the
$W$-orbit of $\tht$, and for each $\eta$, $w'$ ranges over a set of
representatives of $W_{t'}\backslash W /W_t(\eta)$. Summing over the
$W$-conjugacy classes of $T$, we see that in the sum for $\Phi(t,t')$, we can
instead let $\eta$ range over a set of representatives of $\widehat T/ W_t$ and
for each $\eta$, $w'$ range over a set of representatives of
$W_{t'}\backslash W / W_t(\eta)$. Recalling the notation $\Phi(t,t')$
introduced in the proof of Proposition~\ref{prop:t1} we obtain
$$\begin{aligned}
\Phi(t,t')&=\displaystyle\sum_{\tht\in\widehat{T}/W} \Phi_\tht(t,t') \\
  &= a_{t,t'} \displaystyle\sum_{\tht\in\widehat{T}/W} \sum_{\eta \in (\tht^W)/W_t} \sum_{w'\in W_{t'}\backslash W/W_t(\eta)}\ \eta(t ({t'}^{w'})^{-1})\, h(\eta,w') \\
  &= a_{t,t'} \displaystyle\sum_{\eta\in\widehat{T}/W_t} \sum_{w' \in W_{t'}\backslash W/W_t(\eta)}\ \eta(t ({t'}^{w'})^{-1})\, h(\eta,w') \\
  &= a_{t,t'} \displaystyle\sum_{\eta\in\widehat{T}/W_t} \sum_{z\in W/W_t} \sum_{w'\in(W_{t'}\backslash W/W_t(\eta))\cap zW_t}\!\eta(t ({t'}^{w'})^{-1})\, h(\eta,w') \\
  &= a_{t,t'}\! \displaystyle \sum_{z\in W/W_t} \sum_{\eta\in\widehat{T}/W_t} \sum_{w' \in(W_{t'^z}\cap W_t) \backslash W_t /W_t(\eta))}\! \eta(t ({t'}^{zw'})^{-1})\, h(\eta,zw') \\
  &= \displaystyle\sum_{z\in W/W_t} \Phi^{W_t}(t,{t'}^z).
\end{aligned}$$
To see the penultimate equality, let $z\in W$, $\eta\in\widehat T$ and $J$ a set
of representatives of $W_{t'}\backslash W /W_t(\eta)$. Then as $w$ runs
through $J\cap zW_t$, $z^{-1}w$ runs through a set of representatives of
$(W_{{t'}^{z}}\cap W_t) \backslash W_t /W_t(\eta)$. 

Since $t$ is central in $W_t$, the result follows by Proposition~\ref{prop:t1}
and noting that $t'^{z}=t$ for at most one index $z$.
\end{proof}

\begin{rem}   \label{rem:part(2)}
Since $\cH(W,\bu)$ is conjectured to always be strongly symmetric, the
conclusion of Theorem~\ref{thm:B02} is expected to hold without that
assumption. In fact, our proof only required $\cH(W,\bu)$ to satisfy part~(2)
of the definition of strong symmetry in \cite[Def.~3.4]{KMS2}, not~(1), and
this latter statement is related to the parabolic trace conjecture in
\cite{CC20}. 
\end{rem}

\section{Partial character tables for principal blocks of $\ZZ_\ell$-spetses}   \label{sec:spetschars}

In this section we define partial character tables for principal blocks of
$\ZZ_\ell$-spetses in the coprime case, as a special case of the deformed
character tables for $\ZZ_\ell$-reflection groups introduced in the previous
section.

\subsection{The partial character table and block orthogonality}   \label{subsec:chartab spets}
Throughout, $\GG=(W,\La)$ is a simply connected $\ZZ_\ell$-spets such that
$\ell>2$ with $(|W|,\ell)=1$. Note that by \cite[Prop.~2.6]{KMS}, $\ell$ is
very good for $\GG$. Let $q\in1+\ell\ZZ_\ell$ be a prime power and let
$T=\La/\ell^a \La$ be the (abelian) Sylow $\ell$-subgroup of~$\GG(q)$ (see
Lemma~\ref{lem:phi_t}).

Let $B_0=B_0(\GG(q))$ be the principal block of $\GG(q)$, as defined in
\cite[Def.~6.7]{KMS}, \cite[Def.~4.1]{KMS2}. By \cite[Lemma~4.10]{KMS2} there
is a degree-preserving bijection between the set of characters of $B_0$ and the
set of characters of $B_0(\GG_{\ba, a})$ (see Definition~\ref{def:B}) when
$\psi_\ba$ is the admissible \textit{spetsial specialisation} $\psi_\spets$
(see \cite[\S 3B]{KMS2}) for which $a_{rj} =1$ if $j =0$ and $a_{rj} =0$ if
$j >0$, and $\ell^a$ is the highest power of $\ell$ dividing $q-1$. 

Here, as in Section~\ref{sec:chararbitrary} we work under a
fixed $W$-equivariant bijection $\hat{}:T\to\Irr(T)$, $t\mapsto\hat t$.

\begin{defn}   \label{def:chartab spets}
 For $\GG$ as above and $\ga_{\tht,\vhi}\in\Irr(B_0)$ (with
 $\tht\in\Irr(T)$ and $\vhi\in\Irr(W(\tht))$) and $t\in T$, the \emph{character
 value $\ga_{\tht,\vhi}(t)$ of $\ga_{\tht,\vhi}$ at $t$} is given as
 in Definition~\ref{def:chartab}, with $\psi_\ba=\psi_\spets$. We call the
 matrix $\cX=\cX(B_0)=(\ga_{\tht,\vhi}(t))$ the \emph{(partial) character
 table of $B_0$}.
\end{defn} 
Note that when $\psi_\ba=\psi_\spets$, then by \cite[\S4.1]{KMS2}, the degrees
$\ga_{\tht,\la}^{W_{t^w}}(1)$ appearing in Definition~\ref{def:chartab} all
lie in $\ZZ_\ell[x]$ and hence the entries of the character table above all lie
in $\cO[x] $, where $\cO =\ZZ_\ell[\ze]$ for $\ze$ a primitive $\ell^a$th
root of unity.

We anticipate that the analogues of many suitable group- or character-theoretic
statements should hold for this definition. For example, the following asserts
that, regarded as a class function on $T$, $\chi$ specialises to an
$\cO$-linear combination of characters of the underlying Sylow $\ell$-subgroup:

\begin{conj}   \label{conj:restriction}
 For each $\chi\in B_0$, the class function $\chi_{x=q}: T\to\cO$ given by
 $\chi_{x=q}(t):=\chi(t)|_{x=q}$ for $t\in T$ lies in $\cO[\Irr(T)]$.
\end{conj}

We also expect that the analogue of Theorem \ref{t:frob} holds.

\begin{conj}   \label{c:frobspets}
 For each $\chi\in\Irr(B_0)$ we have
 $$\Big(\sum_{t\in T/W} |\GG:C_{\GG}(t)|\chi(t)\Big)_{x=q}\equiv 0\pmod{|T|}.$$
\end{conj}

Conjectures \ref{conj:restriction} and \ref{c:frobspets} both hold when $W$ is
rational as will be seen in Proposition \ref{p:blockbijection} below and are
easily checked for small non-rational examples, such as when $W$ is cyclic of
order~$e >2$. We defer a thorough investigation to future work. Here, our main
conjecture of interest is the following analogue, for spetses, of block
orthogonality for the principal $\ell$-block of a finite group
(Theorem~\ref{t:centblock}):

\begin{conj}   \label{cor:main1}
 For all $t,t'\in T$,
 $$\sum_{\chi \in \Irr(B_0)} \chi(t){\chi({t'}^{-1})}
    = \begin{cases} \dim(B_0(t))& \text{if $t\sim_W t'$,}\\
                    0& \text{otherwise,}\end{cases}$$
 where $B_0(t)$ denotes the principal block of $C_\GG(t)$.
\end{conj}

We can now prove:

\begin{thm}   \label{thm:blockorth}
 Conjecture~\ref{cor:main1} holds if $\cH(W,\bu)$ is strongly symmetric.
\end{thm}

\begin{proof}
By definition, for any non-conjugate $t,t'\in T$ the considered scalar product
in Conjecture~\ref{cor:main1} equals the corresponding product
for~$\GG_{\ba,a}$ in Theorem~\ref{thm:B02} and thus vanishes, if $\cH(W,\bu)$
is strongly symmetric.
\end{proof}

\subsection{The case of finite reductive groups}   \label{sec:red group}
In this section we show that the character values introduced in
Definition~\ref{def:chartab spets} specialise to those of actual finite groups
of Lie type when the underlying Weyl group is rational. Let $\bG$ be a connected
reductive linear algebraic group of simply connected type over an algebraically
closed field of characteristic~$p$ and $F:\bG\to\bG$ a split Frobenius map with
respect to an $\FF_q$-structure, for $q$ a power of $p$. Let $\bT$ be a
maximally $F$-split torus of~$\bG$ and let $W=N_\bG(\bT)/\bT$, the Weyl group
of $\bG$. Let $T$ be the Sylow $\ell$-subgroup of $\bT^F$ and let $C$ be the
principal $\ell$-block of $\bG^F$. Suppose that $\ell$ divides $q-1$ and does
not divide $|W|$. Then $T$ is a Sylow $\ell$-subgroup of $\bG^F$. 
Let $\GG=(W,\La)$ be the corresponding $\ZZ_\ell$-spets and $B_0$ be the
principal block of $\GG(q)$ where $\psi_\ba$ is the spetsial
specialisation and $\ell^a=(q-1)_\ell$. Note that $T$ may be identified with
the Sylow $\ell$-subgroup of $\GG(q)$ (see \cite[\S 5.3]{KMS}).

\begin{prop}   \label{p:blockbijection}
 There exists a bijection $\Theta:\Irr(B_0)\to\Irr (C)$ for which
 $\ga(t)|_{x=q} = \Theta(\ga)(t)$ for all $\ga\in\Irr(B_0)$ and $t\in T$.
\end{prop}

\begin{proof}
For $\tht\in\Irr(\bT^F)$ and $\bL$ a $1$-split Levi subgroup of $\bG$
containing $\bT$ denote by $\cE(\bL^F,(\bT,\tht))$ the 1-Harish-Chandra series
of $\bL^F$ above $(\bT,\tht)$ and let $W_\bL = N_\bL(\bT)/\bT \leq W$ denote
the Weyl group of $\bL$ (so, $W_{\bG}= W$).
By the Howlett--Lehrer comparison theorem \cite[Thm~3.27]{GM20} for a fixed
$\tht\in\Irr(\bT^F)$ there is a collection of bijective labellings
$\psi\mapsto\eta_\psi$ of $\cE(\bL^F,(\bT,\tht))$ by $\Irr(W_{\bL}(\tht))$
where $\bL $ runs through the $1$-split Levi subgroups of $\bG$ containing
$\bT$, such that for any $\bL$ and any $\psi\in\Irr(W_{\bL}(\tht))$,
$$R_{\bL}^\bG(\eta_\psi)=\eta_{\Ind_{W_{\bL}(\tht)}^{W(\tht)}(\psi)},$$
where we have extended $\psi\mapsto\eta_\psi$ linearly.

Restriction induces a $W$-equivariant bijection between the set of $\ell$-power
order elements of $\Irr(\bT^F)$ and $\Irr(T)$. We henceforth identify the two
sets via this bijection and let $I$ be a set of representatives of the
$W$-orbits of $\Irr(T)$. 
By \cite[Sec.~3]{BrMi}, $\Irr(C)$ is the disjoint union of the 1-Harish-Chandra
series $\cE(\bG^F,(\bT,\tht))$, where $\tht\in\Irr(\bT^F)$ runs over a set of
representatives of the $W$-classes of the characters of $\ell$-power order.
Thus, $\Irr(C)$ is the disjoint union of the series $\cE(\bG^F,(\bT,\tht))$,
$\tht\in I$. For $\ga =\ga_{\tht,\vhi}\in\Irr(B_0)$ with $\tht\in I$ and
$\vhi\in \Irr(W(\tht))$, let $\Theta(\ga_{\tht,\vhi}):=\eta_\vhi$ be the
element of $\cE(\bG^F,(\bT, \tht))$ labelled by $\vhi$. Then
$\Theta:\Irr (B_0)\to\Irr(C)$ is a bijection.

Now suppose $\tht\in I$, $\vhi\in\Irr(W(\tht))$ and let
$\chi :=\Theta (\ga_{\tht,\vhi}) =\eta_\vhi$. Let $t\in T$,
$\bL:=\bC_\bG(t)$, a $1$-split Levi subgroup of $\bG$ containing $\bT$, and let
$W_t := W_{\bL}$. By \cite[Cor.~3.3.13]{GM20} we have
$$\chi(t)=\sum_{\la\in\Irr(\bL^F)}\bla\chi,\RLG(\la)\bra\,\la(t).$$
Assume that $\la\in\Irr(\bL^F)$ contributes to the above sum, that is,
$\langle\chi,\RLG(\la)\rangle\ne0$. Then, by Harish-Chandra theory, $\la$ lies
in the 1-Harish-Chandra series of $(\bT,\tht)^g$ for some $g\in\bG^F$, and we
have
$$\chi(t)=\sum_{X\in\cM}\sum_{\la\in\cE(\bL^F,X)}
  \bla\chi,\RLG(\la)\bra\,\la(t)$$
where $\cM:=\{(\bT^g,\tht^g)\mid g\in\bG^F,\,\bT^g\le\bL\}/\sim_{\bL^F}$ and
where $\cE(\bL^F,X)$ for $X\in\cM$ denotes the 1-Harish-Chandra series of
$\bL^F$ above $X$. Now, $\cM$ is in bijection with $W_t\backslash W/W(\tht)$,
so
$$\begin{aligned}
\chi(t)
  &=\sum_{w\in W_t\backslash W/W(\tht)}\sum_{\la\in\cE(\bL^F,\tw{w}(\bT,\tht))}
   \bla\chi,\RLG(\la)\bra\,\la(t)\\
  &=\sum_{w\in W_t\backslash W/W(\tht)}\sum_{\la\in\cE((\bL^w)^F,(\bT,\tht))}
   \bla\chi,R_{\bL^w}^\bG(\la)\bra\,\la(t^w).
\end{aligned}$$
Applying the comparison theorem we obtain
$$\chi(t)=\sum_{w\in W_t\backslash W/W(\tht)}\sum_{\psi\in\Irr(W_{t^w}(\tht))}
 \bla\vhi,\Ind_{W_{t^w}(\tht)}^{W(\tht)}(\psi)\bra\,\eta_\psi(t^w),$$
where $W_{t^w}:= W_{\bL^w}= W_t^w$.

Now $t^w$ is central in $\bL^w$ and so for all
$\eta\in\cE(\bL^{wF},(\bT,\tht))$ we
have $\eta(t^w)=\tht(t^w)\eta(1)$ (see \cite[Lemma~2.2]{Ma07}).
Combining with the above we thus obtain
$$\chi(t)
  =\sum_{w\in W_t\backslash W/W(\tht)}\tht(t^w)\sum_{\psi\in\Irr(W_{t^w}(\tht))}
   \bla\vhi,\Ind_{W_{t^w}(\tht)}^{W(\tht)}(\psi)\bra\,\eta_\psi(1).$$
Now by \cite[Prop.~(68.30)(iii)]{CR81} the Schur element of $\cH(W,\bu)$
corresponding to $1_W$, specialised at $x=q$, is $|\bG^F:\bT^F|_{p'}$ (see also
\cite[Prop.~4.6 and Lemma~4.10]{KMS2}). This, in combination with
\cite[Def.~3.2.13, Thm~3.2.18]{GM20}, shows that
$\eta_\psi(1)= \ga_{\tht,\psi}^{W_{t^w}}(1)_{x=q}$ for $\ga_{\tht,\psi}$ as in
Definition \ref{def:B}. We thus obtain the formula in
Definition~\ref{def:chartab}.
\end{proof}

We end this subsection with the following inductive reformulation of Definition
\ref{def:chartab spets} in terms of 1-Harish-Chandra theory, see
\cite[Ax.~4.16]{BMM14} for this notion in the case of spetses. The argument is
essentially a reversal of that given in the second part of the proof of
Proposition~\ref{p:blockbijection}.

\begin{lem}   \label{l:conncharform}
 Suppose that $t\in T$ and $\chi\in\Irr(B_0)$. Then $\LL:=C_\GG(t)$ is a
 $1$-split Levi of $\GG$. Assume $\chi$ lies in the $1$-Harish-Chandra series
 above the unipotent $1$-cuspidal pair $(\TT,\tht)$ for $\TT$ a Sylow
 $1$-torus. Let
 $\cM:=\{^w(\TT,\tht)\mid w\in W \}/W_\LL$. Then
 $$\chi(t)=\sum_{\la\in\cE(\LL,X),\,X\in\cM}
           \langle\rho,R_\LL^\GG(\la)\rangle\,\la(t).$$
           Moreover, the value
 $\big(|\GG:C_\GG(t)|\chi(t)/\chi(1)\big)_{x=q}$ is integral over $\ZZ_\ell$.
\end{lem}

\begin{proof}
Since $\ell$ is very good for $\GG$ and $\ell\nmid |W|$, $\LL= (C_W(\La'),\La)$
where $\La'= \La^{C_W(t)}$ (see the proof of \cite[Prop.~2.3]{KMS}). In
particular, $\LL$ is a Levi of $\GG$. The first two lines of the
proof of
\cite[Lemma 4.10]{KMS2} argue that since $q\equiv1\pmod\ell$, $\LL$ contains a
Sylow 1-torus. Since $\phi=1$, a Sylow 1-torus of $\GG$ is a maximal torus.
So $Z(\LL)$ is a 1-torus, and $\LL$, being its centraliser, is $1$-split.
Assume that $\chi=\ga_{\tht,\vhi}$ with $\vhi \in \Irr(W(\tht))$.
Since $\cM$ is in bijection with $W_t\backslash W/W(\tht)$, we have
$$\begin{aligned}
\sum_{X\in\cM}\sum_{\la\in\cE(\LL,X)}
  \bla\chi,R_\LL^\GG(\la)\bra\,\la(t)
  &=\sum_{w\in W_t\backslash W/W(\tht)}\sum_{\la\in\cE(\LL,\tw{w}(\TT,\tht))}
   \bla\chi,R_\LL^\GG(\la)\bra\,\la(t)\\
  &=\sum_{w\in W_t\backslash W/W(\tht)}\sum_{\la\in\cE(\LL^w,(\TT,\tht))}
   \bla\chi,R_{\LL^w}^\GG(\la)\bra\,\la(t^w)
    \end{aligned}$$
where $W_{t^w}:= W_{\LL^w}= W_t^w$. Now $t^w$ is central in $\LL^w$ and so for
$\la\in\cE(\LL^w,(\TT,\tht))$ we have $\la(t^w)=\tht(t^w)\la(1)$ (see
Example~\ref{ex:central}). Moreover, $\lambda=\ga_{\tht,\psi}^{W_{t^w}}$
for some $\psi \in \Irr(W_{t^w}(\tht))$. 
Combining these observations with the above we obtain
$$\begin{aligned}
\sum_{X\in\cM}\sum_{\la\in\cE(\LL,X)}
  \bla\chi,R_\LL^\GG(\la)\bra\,\la(t)
  &=\sum_{w\in W_t\backslash W/W(\tht)}\sum_{\psi\in\Irr(W_{t^w}(\tht))}
 \bla\vhi,\Ind_{W_{t^w}(\tht)}^{W(\tht)}(\psi)\bra\,\ga_{\tht,\psi}^{W_{t^w}}(t^w)\\
  &=  \sum_{w\in W_t\backslash W/W(\tht)}\tht(t^w)\sum_{\psi\in\Irr(W_{t^w}(\tht))}
   \bla\vhi,\Ind_{W_{t^w}(\tht)}^{W(\tht)}(\psi)\bra\,\ga_{\tht,\psi}^{W_{t^w}}(1),
\end{aligned}$$
as in Definition~\ref{def:chartab}.

For the final statement, observe that since $\LL=C_\GG(t)$ is $1$-split,
$|\GG:C_\GG(t)|_{x=q}$ is prime to~$\ell$.
Let $\chi:=\ga_{\tht,\vhi}^{W(\tht)}\in\cE(\GG,(\TT,\tht))$ and assume
$\tht=\widehat{s}$ for some $s \in S$. For some
$w\in W_t \backslash W/ W(\tht)$, let
$\ga_{\tht,\psi}^{W_{t^w}}\in\cE(\LL^w,(\TT,\tht))$
be one of the characters appearing in the above formula for
$\chi(t)$, so $\ga_{\tht,\psi}^{W_{t^w}}$ and $\chi$ have the same
$1$-Harish-Chandra vertex $(\TT,\tht)$. Now,
$$\chi(1)=|\GG:C_\GG(s)|_{x'}\,\ga_{1,\vhi}^{W(\tht)}(1)\quad\text{and}\quad
  \ga_{\tht,\psi}^{W_{t^w}}(1)
  =|\LL^w:C_{\LL^w}(s)|_{x'}\,\ga_{1,\psi}^{W_{t^w}(\tht)}(1).$$
Thus $|\GG:C_\GG(t)|\chi(t)/\chi(1)$ is equal to

$$ \frac{|C_\GG(s)|_{x'}}{|\LL|_{x'}\,\ga_{1,\vhi}(1)}\sum_{w\in W_t\backslash W/W(\tht)}\tht(t^w)\sum_{\psi\in\Irr(W_{t^w}(\tht))}
   \bla\vhi,\Ind_{W_{t^w}(\tht)}^{W(\tht)}(\psi)\bra\,|\LL^w:C_{\LL^w}(s)|_{x'}\,\ga_{1,\psi}^{W_{t^w}(\tht)}(1),$$
so that on specialisation at $x=q$ it suffices to see that
$\ga_{1,\psi}^{W_{t^w}(\tht)}(1)_{x=q}$ and
$\ga_{1,\vhi}^{W(\tht)}(1)_{x=q}$
have the same $\ell$-adic valuations. But by $1$-Harish-Chandra theory,
both have the same $\ell$-adic valuation as $\tht(1)_{x=q}=1$, so the claim
follows. 
\end{proof}

\subsection{Steinberg characters}   \label{sec:steinberg}
To conclude this section, we investigate the character values in a particular
case. As in \S\ref{subsec:chartab spets}, we assume here that $\ell> 2$ is
relatively prime to $|W|$ and $q\in 1+\ell\ZZ_\ell$.

\begin{lem}   \label{lem:val linear}
 Let $\tht=1$ and $\vhi\in\Irr(W)$ be linear. Then the values of the unipotent
 character $\ga_{1,\vhi}$ are given by
 $$\ga_{1,\vhi}^W(t)=\ga_{1,\vhi|_{W_t}}^{W_t}(1)\qquad
   \text{ for all $t\in T$}.$$
\end{lem}

\begin{proof}
In our situation, $W(s)=W$, so the first sum in the definition of
$\ga_{1,\vhi}^W(t)$ has just one term corresponding to $w=1$, and $\widehat s$
is the trivial character. Furthermore, for $\la\in\Irr(W_t)$
$$\bla\vhi,\Ind_{W_t}^{W}(\la)\bra
  =\bla\Res_{W_t}^W(\vhi),\la\bra=\begin{cases}
    1& \text{if $\la=\vhi|_{W_t}$,}\\ 0& \text{else,}\end{cases}$$
since the restriction of the linear character $\vhi$ to any subgroup is
irreducible. The stated formula follows.
\end{proof}

This gives an analogue, or rather, generalisation, of the well-known formula
for the values of the Steinberg character of groups of Lie type:

\begin{cor}   \label{cor:St}
 Assume $W$ is spetsial and generated by reflections of order~$2$ and let
 $\eps\in\Irr(W)$ be the sign character (taking value~$-1$ on all
 reflections). Let $\psi_\ba$ be the spetsial specialisation. Then
 $$\ga_{1,\eps}^W(t)= x^{N(W_t)}\text{ for all $t\in T$},$$
 where $N(W_t)$ is the number of reflections in $W_t$.
\end{cor}

\begin{proof}
By \cite[6B]{Ma95} and the tables in \cite{BMM14}, in our situation we have
that $\ga_{1,\eps}^Y(1)=x^{N(Y)}$ for any irreducible spetsial reflection
group~$Y$.
(This is an instance of Alvis--Curtis duality on unipotent character degrees.)
Since parabolic subgroups of spetsial groups are again spetsial, and any
spetsial group is a direct product of irreducible ones, and
as $N(Y_1\times Y_2)=N(Y_1)+N(Y_2)$, the claim is immediate from
Lemma~\ref{lem:val linear}.
\end{proof}

The following example provides further evidence for
Conjecture~\ref{conj:restriction} when $\St:=\ga_{1,\eps}^W$ is the Steinberg
character. 

\begin{exmp}
 Let $W=G(e,e,3)$, with $e\ge3$, an imprimitive spetsial complex reflection
 group generated by reflections of order 2. Assume $2<\ell\equiv1\pmod e$, so
 $W$ is a $\ZZ_\ell$-reflection group, and let $q\in 1+\ell^a\ZZ_\ell^\times$.
 There are five types of parabolic subgroup $W_0$ of~$W$: apart from $W$,
 $I_{2e}$ occurs three times as centraliser, $A_2$ occurs $e^2$ times, $A_1$
 occurs $3e$ times and $A_0$ once.  The following table records the number of
 $t\in T$ with $W_t=W_0$ together with the value $\St(t)$ for such elements:
{\small $$\begin{array}{c|c|c|c|c|c}
 W_0 & W & I_{2e} & A_2 & A_1 & A_0 \\ \hline
 |\{t \in T: W_t=W_0\}| & 1 & \ell^a-1 & \ell^a-1 & (\ell^a-1)(\ell^a-e) & (\ell^a-1)(\ell^a-2(e-1))(\ell^a-e-1) \\
 \St(t) & x^{3e} & x^e & x^3 & x & 1
\end{array}$$}
From this information we obtain that $|T|\,\langle\St,1_T\rangle$ equals
$\ell^{3a}+3e(x-1)\ell^{2a}+p_1\ell^a+p_2$ with
$$p_1:=3(x^e-1)-3(x-1)e+e^2(x-1)^2(x+2),\quad
  p_2:=x^{3e}-1-3e(x-1)-p_1.$$
Note that $p_1$ has a double zero at $x=1$, and $p_2$ has a triple zero
(differentiate once or twice), so indeed we get
$|T|\,\langle\St,1_T\rangle|_{x=q}\equiv0\pmod {\ell^{3a}}$ as predicted by
Conjecture~\ref{conj:restriction}.
\par
Now suppose that $\widehat s\in\Irr(T)$ is regular so that $\widehat s$
restricts non-trivially to each of the $e^2+3$ one-dimensional subgroups of $T$
fixed by parabolic subgroups of rank~2. Then again
$$|T|\,\langle\St,\widehat s\rangle|_{x=q}=p_2(q)\equiv 0\pmod{\ell^{3a}}.$$
The other possibilities for $\widehat{s}$ may be handled similarly.
\end{exmp}

\newpage
\section{Values of unipotent characters}   \label{sec:valuesunip}

\subsection{Deligne--Lusztig characters and unipotent characters}
Let $\GG=(W\phi,\La)$ be a $\ZZ_\ell$-spets.  Let $\cO$ be an integral extension
of $\ZZ_\ell$ containing a primitive $|\phi|$-th root of unity, with field of
fractions $K$. By \cite{Ma95} and \cite[Thm~8.4]{Ma00} in conjunction with
\cite[Ax.~4.1]{BMM14} associated to a ($\ZZ_\ell$-) spets $\GG=(W\phi,\La)$
there comes a set of \emph{unipotent characters $\Uch(\GG)$}, together with a
\emph{degree} function
$$\Uch(\GG)\longrightarrow\QQ_\ell[x],\quad \rho\mapsto \rho(1),$$
and an action of $N_{\GL(\La)}(W\phi)$ on $\Uch(\GG)$ (see also \cite{BMM14}
for explicit tables for primitive~$W$). Here we note that for any
$\ga\in\Uch(\GG)$, $\ga(1)$ is a polynomial with coefficients in the character
field of $\langle W,\phi\rangle$, and consequently, since $\GG$ is a
$\ZZ_\ell$-spets, $\ga(1)\in\QQ_\ell[x]$. If $W$ is a Weyl group, $\Uch(\GG)$
is just a labelling set for the unipotent
characters of an associated finite reductive group of type $W\phi^{-1}$ defined
over $\FF_q$ together with their degree polynomials. Let $\cU(\GG)$ be the
$K$-vector space with orthonormal basis $\Uch(\GG)$. By linear
extension to any $\ga\in\cU(\GG)$ is attached its degree $\ga(1)\in K[x]$.

For $\vhi\in\Irr(W)^\phi$, a $\phi$-invariant irreducible character of~$W$,
we choose an extension $\tvhi$ to $W\langle\phi\rangle$ and let $\Irr(W\phi)$
denote the set of restrictions of these $\tvhi$ to the coset~$W\phi$. (These
are hence defined up to multiplication by $o(\phi)$th roots of unity.)
Then the space $\cU(\GG)$ possesses a subspace $\cU_1(\GG)$ with an
orthonormal basis consisting of \emph{uniform almost characters} $R_\vhi$,
$\vhi\in\Irr(W\phi)$, which are defined as certain explicit linear combinations
of the unipotent characters via a Fourier transform matrix (see
\cite[4A and~6C]{Ma95} for the imprimitive groups, and \cite{MChev} for
the primitive ones). From these we define certain elements $R_w^\GG$,
$w\in W$, called \emph{Deligne--Lusztig characters}, by
$$R_w^\GG:=\sum_{\vhi\in\Irr(W\phi)}\vhi(w\phi) R_\vhi.$$
Note that this depens only on the $\phi$-conjugacy class $[w]_{W\phi}$ of $w$
Thus, by orthogonality of $\phi$-class functions we have
$$R_\vhi=|W|^{-1}\sum_{w\in W} \vhi(w\phi)R_w^\GG\quad\in\cU(\GG),$$
and the orthonormality of the $R_\vhi$ implies that the $R_w$, for $w\phi$
running over a system of representatives of the $W$-classes in~$W\phi$, form an
orthogonal basis of $\cU_1(\GG)$ satisfying
$\langle R_v^\GG,R_w^\GG\rangle=\delta_{[v],[w]}|C_{W\phi}(w)|$ for all
$v,w\in W$. We call $\cU_1(\GG)$ the subspace of $\cU(\GG)$ of \emph{uniform
functions}, and let $\cU_2(\GG)$ denote its orthogonal complement in $\cU(\GG)$.

Assume $S$ is a set equipped with a map $t\mapsto W(t)w_t$ which to any $t\in S$
associates a $\phi$-stable coset $W(t)w_t$ of a subgroup $W(t)$ in $W$ such
that $W(t)^\circ$ is spetsial. (Recall from Definition~\ref{def:2.5} that
${}^\circ$ means taking the subgroup generated by reflections.) We then
set $C_\GG(t):=(W(t)w_t\phi,\La)$. We propose to define values of the unipotent
characters of $\GG$ on~$S$, by defining the values of the Deligne--Lusztig
characters $R_w^\GG$ for $w\in W$:

\begin{defn}   \label{def:val DL}
 Let $t\in S$ and $\LL=(W_\LL w_\LL\phi,\La):=C_\GG(t)^\circ$. Then
 $$R_w^\GG(t):=|C_{W\phi}(w)|\cdot|W_\LL|^{-1}\sum_v R_v^\LL(1)\ 
   \in K[x]$$
 where the sum runs over elements $v\in [w]_{W\phi}\cap W_\LL w_\LL$, that
 is, elements of $W_\LL w_\LL$ in the $\phi$-conjugacy class of $w$.
 We decree that all $\rho\in\cU_2(\GG)$ take value~0 on $t$.
 By linear extension this defines values of any $\rho\in\cU(\GG)$ on $S$.
\end{defn}

Observe that this definition only depends on $C_\GG(t)^\circ$, not on $t$
itself. Also, the spets $\LL$ is not necessarily a Levi sub-reflection datum
of~$\GG$, but just a sub-reflection datum in the sense of
\cite[Def. 3.1]{BMM99}. By grouping together the summands for conjugate
elements $v$, the above formula can be rewritten
$$R_w^\GG(t)=\frac{|C_{W\phi}(w)|\cdot
  |[w]_{W\phi}\cap W_\LL w_\LL|}{|W_\LL|}R_v^\LL(1).$$

This definition is of course motivated by the following:

\begin{exmp}   \label{exmp:DL}
 Assume $\GG$ is a rational spets and $(\bG,F)$ is a corresponding finite
 reductive group over a field of characteristic~$p$. Let $S\subset G:=\bG^F$
 be the set of semisimple elements. For $t\in S$, its centraliser is again
 reductive and we let $W(t)w_t$ denote the corresponding reflection coset
 (which determines the rational type of $C_\bG^\circ(t)$). Then (see
 \cite[Prop.~7.5.3]{Ca}) the values of Deligne--Lusztig characters of $G$ are
 given by
 $$R_w^\bG(t)=\eps_w\eps_\bL|T_w|^{-1}|\bL^F|_p^{-1}\sum_g 1
   \qquad\text{for any semisimple $t\in G$},$$
 where $\bL:=C_\bG^\circ(t)$ and the sum runs over elements $g\in G$ with
 $g^{-1}tg\in T_w:=\bT_w^F$, that is, with $t\in gT_wg^{-1}$, that is, with
 $g\bT_wg^{-1}\le \bL$. Let us count such $g$. For this, we count suitable
 $F$-stable maximal tori $\bT$ of $\bL$. Any such $\bT$ is stabilised by
 $|N_G(\bT)|$ elements $g\in G$. Furthermore, $\bT$ has $|L:N_L(\bT)|$ many
 distinct conjugates under $L:=\bL^F$. So any class of $F$-stable maximal tori
 in $\bL$, indexed
 by $v\in W_\bL$ say, which contains one of the $g\bT_wg^{-1}$ will account for
 $|N_G(\bT)|\cdot|L:N_L(\bT)|=|L|\cdot|C_{WF}(v):C_{W_\bL F}(v)|$ elements~$g$
 (using that $|N_L(\bT_v)|=|T_v|\cdot|C_{W_\bL F}(v)|$). The relevant $v$ are
 those lying in $[w]\cap W_\bL$, and any such~$v$ lies in a $W_\bL F$-class of
 length $|W_\bL:C_{W_\bL F}(v)|$. Furthermore, note that
 $\eps_w\eps_\bL|L|_{p'}/|T_v|=R_v^\bL(1)$. Piecing all this
 together we find that 
 $$R_w^\bG(t)=|C_{WF}(w)|\cdot|W_\bL|^{-1}\sum_v R_v^\bL(1),$$
 as in Definition~\ref{def:val DL}.   \par
 Moreover, by a basic result in Deligne--Lusztig theory \cite[Prop.~7.5.5]{Ca},
 the characteristic function of any semisimple class of $G$ is uniform, so any
 class function orthogonal to all uniform functions vanishes on $t$. Thus,
 Definition~\ref{def:val DL} is a generic version of the formula for rational
 spetses.
\end{exmp}

Returning to the general situation, we now derive a formula for the value of
arbitrary $\rho\in\cU(\GG)$, for which we use (compare to
\cite[Prop.~10.2.4]{DiMi2}):

\begin{lem}   \label{lem:uniform proj}
 Let $\rho\in\cU(\GG)$. The orthogonal projection $\rho_\uni$ of $\rho$ onto
 $\cU_1(\GG)$ (its \emph{uniform projection}) is given by
 $$\rho_\uni=|W|^{-1}\sum_{w\in W}\langle\rho,R_w^\GG\rangle\,R_w^\GG.$$
\end{lem}

\begin{proof}
Write $\rho_\uni=\sum_{w\in W}a_w R_w^\GG$, with $a_w\in K$ constant on
$\phi$-classes. Then for any $v\in W$ we have
$$\bla\rho,R_v^\GG\bra=\bla\rho_\uni,R_v^\GG\bra
  =\sum_{w\in W}a_w\bla R_w^\GG,R_v^\GG\bra.$$
Using that $\langle R_w^\GG,R_v^\GG\rangle=\delta_{[v],[w]}|C_{W\phi}(v)|$
(see above) this yields $\langle\rho,R_v^\GG\rangle=|W|a_v$, so the
multiplicity $a_w$ is as claimed.
\end{proof}

\begin{cor}   \label{cor:val unip}
 Let $\rho\in\cU(\GG)$ and $t\in S$ with
 $\LL=(W_\LL w_\LL\phi,\La):=C_\GG(t)^\circ$. Then
 $$\rho(t)
   =|W_\LL|^{-1}\sum_{v\in W_\LL}\bla\rho,R_v^\GG\bra\, R_v^\LL(1).$$
\end{cor}

\begin{proof}
By Definition~\ref{def:val DL}, $\rho(t)=\rho_\uni(t)$, so by
Lemma~\ref{lem:uniform proj},
$$\rho(t)=|W|^{-1}\sum_{w\in W}\langle\rho,R_w^\GG\rangle\,
  |C_{W\phi}(w)|\cdot|W_\LL|^{-1}\sum_v R_v^\LL(1),$$
with $v$ running over $[w]_{W\phi}\cap W_\LL w_\LL$. We replace the first sum
by a sum over $\phi$-conjugacy classes (any of which contains
$|W:C_{W\phi}(w)|$ elements) to get
$$\rho(t)=\sum_{[w]\subset W}\bla\rho,R_w^\GG\bra\,
  |W_\LL|^{-1}\sum_{v} R_v^\LL(1)
  =|W_\LL|^{-1}\sum_{v\in W_\LL}\bla\rho,R_v^\GG\bra\,R_v^\LL(1).
  \qedhere$$
\end{proof}

Let us recall some definitions and facts from \cite{BMM14} and \cite{Ma26}.
For any root of unity $\ze\in\ZZ_\ell^\times$ and $\Phi:=x-\ze$, there is a
subset of \emph{$\Phi$-cuspidal} elements of $\Uch(\GG)$. The $\Phi$-cuspidal
$\rho\in\Uch(\GG)$ are characterised by the property that 
\begin{equation}\label{eq:cusp}
  \rho(1)_\Phi=\frac{|\GG|_\Phi}{|Z^\circ(\GG)|_\Phi}
  =|\GG_\der|_\Phi\end{equation}
(see \cite[Def.~4.5]{BMM14}) where 
$$\GG_\der:= \big(W\phi,([\Lambda,W]\otimes_{\ZZ_\ell}\QQ_\ell)\cap \Lambda\big), \qquad [\Lambda,W]:=\langle (w-1)\Lambda \mid w \in W \rangle.$$
A \emph{$\Phi$-cuspidal pair} of $\GG$ is a pair $(\LL,\la)$ where $\LL$ is a
$\Phi$-split Levi of $\GG$ and $\la\in\Uch(\LL)$ is $\Phi$-cuspidal. 
For any root of unity $\ze\in\ZZ_\ell^\times$, by \cite[Thm~2.1]{Ma26},
$\Uch(\GG)$ is partitioned into $\Phi$-Harish-Chandra series above
$\Phi$-cuspidal pairs $(\LL,\la)$, in bijection with $\Irr(W_\GG(\LL,\la))$
and with degrees given in terms of Schur elements of the Hecke algebra
$\cH_\GG(\LL,\la)$ associated to $W_\GG(\LL,\la)$. For $(\LL,\la)$ a
$\Phi$-cuspidal pair of $\GG$, we denote by
$\cE(\GG,(\LL,\la))\subseteq \Uch(\GG)$ the corresponding $\Phi$-Harish-Chandra
series. The formula in Corollary~\ref{cor:val unip} can then be reformulated as
follows, assuming that $\LL$ is $\Phi$-split for some~$\Phi$:

\begin{cor}   \label{cor:Curtis}
 Let $\rho\in\cU(\GG)$ and $t\in S$ with $\LL:=C_\GG(t)^\circ$. Assume $\LL$ is
 a $\Phi$-split Levi of~$\GG$ for some $\Phi$, and $\rho$ lies in the
 $\Phi$-Harish-Chandra series above the unipotent $\Phi$-cuspidal pair
 $(\MM,\mu)$. Let $\cM:=\{(\MM,\mu)^w\mid w\in W,\ \MM^w\le\LL\}/W_\LL$. Then
 $$\rho(t)=\sum_{\la\in\cE(\LL,X),\,X\in\cM}
           \langle\rho,R_\LL^\GG(\la)\rangle\,\la(1).$$
\end{cor}

\begin{proof}
By Corollary~\ref{cor:val unip} and Lemma~\ref{lem:uniform proj} we have
$$\begin{aligned}
  \rho(t)=|W_\LL|^{-1}\sum_{w\in W_\LL}\bla\rho,R_w^\GG\bra\,R_w^\LL(1)
  =&|W_\LL|^{-1}\sum_{w\in W_\LL}\bla\tw{*}R_\LL^\GG(\rho),R_w^\LL\bra\,R_w^\LL(1)\\
  =&\tw{*}R_\LL^\GG(\rho)_\uni(1)=\tw{*}R_\LL^\GG(\rho)(1).
\end{aligned}$$
Here, $\tw{*}R_\LL^\GG$ is the adjoint of $R_\LL^\GG$ with respect to our
scalar product. Now by $\Phi$-Harish-Chandra theory \cite[Thm~2.1(b)]{Ma26}
$$\tw{*}R_\LL^\GG(\rho)=\sum_{\la\in\cE(\LL,X),\,X\in\cM}
\langle\tw{*}R_\LL^\GG(\rho),\la\rangle\,\la=\sum_{\la\in\cE(\LL,X),\,X\in\cM}
\langle\rho,R_\LL^\GG(\la)\rangle\,\la,$$
which gives the stated formula.
\end{proof}

The degree $f_\vhi^\GG:=R_\vhi(1)$ of the almost character $R_\vhi$ is called
the \emph{fake degree} of $\vhi\in\Irr(W\phi)$.
By \cite[\S1.2.2]{BMM14} this is the graded multiplicity of $\tvhi$ in the
coinvariant algebra of~$W$. We extend both $R_\vhi$ and $f_\vhi^\GG$ linearly
to $\vhi\in\ZZ\Irr(W\phi)$.

\begin{lem}   \label{lem:val almost}
 Let $\vhi\in\Irr(W\phi)$. Then the almost character $R_\vhi$ takes
 values
 $$R_\vhi(t)=f_\vhi^\LL\qquad\text{for $t\in S$ with
   $\LL=(W_\LL\phi_\LL,\La)=C_\GG(t)^\circ$},$$
 where $f_\vhi^\LL$ is the fake degree of $\vhi|_{W_\LL}$ for $\LL$ (extended
 linearly to $\ZZ\Irr(W_\LL\phi_\LL)$).
\end{lem}

\begin{proof}
By definition we have (arguing as in the proof of Corollary~\ref{cor:val unip})
$$\begin{aligned}
  R_\vhi(t)=|W|^{-1}\sum_{w\in W} \vhi(w\phi)R_w^\GG(t)
  =&\sum_{[w]\subset W}|W_\LL|^{-1}|\,[w]\cap W_\LL w_\LL|\,\vhi(w\phi)\,R_w^\LL(1)\\
  =&|W_\LL|^{-1}\sum_{v\in W_\LL}\vhi(v\phi_\LL)\,R_v^\LL(1)
\end{aligned}$$
which, by definition, is $f_\vhi^\LL$.
\end{proof}

\subsection{Specialisation to the Sylow subgroup}   \label{sec:res_sylow}
We assume $\GG=(W\phi,\La)$ is simply connected and $\ell$ is odd. Let
$q\in\ZZ_\ell^\times$ and let $S$ be the Sylow $\ell$-subgroup of $(\GG,q)$.
As in $\S$\ref{subsec:centr} let $\cF$ be the saturated fusion system of
$(\GG,q)$ on $S$ and recall from \cite[\S 6.2]{KMS} that to each fully
$\cF$-centralised $t\in S$ one may associate its centraliser
$C_\GG(t)=(W(t)\phi_t,\La)$, a spets by \cite[Prop.~6.2]{KMS}. Now let
$s\in S$ and $t$ be a fully $\cF$-centralised $\cF$-conjugate of $s$.
For $\chi\in\cU(\GG)$ define:
\begin{equation}   \label{e:uniformula}
  \chi(s):=\chi(t)=|W_\LL|^{-1}\sum_{v\in W_\LL}\bla\chi,R_v^\GG\bra\, R_v^\LL(1)
\end{equation}
to be the character value given in the conclusion of
Corollary~\ref{cor:val unip}, where $R_w^\GG(t)$ is defined with
respect to the map $t\mapsto W(t)\phi_t$ as in Definition \ref{def:val DL}.

\begin{rem} We make the following observations concerning the formula \eqref{e:uniformula}.
\par
(1) By \cite[Prop. 5.5]{KMS}, $\chi(s)$ is independent of the choice of fully
 $\cF$-centralised $\cF$-conjugate $t$. In particular, $\chi(s)=\chi(s')$
 whenever $s, s' \in S$ are $\cF$-conjugate.
\par
(2) If $\GG$ is generated by true reflections, the Steinberg character
 $\St:=\ga_{1,\eps}^W$ from \S \ref{sec:steinberg} agrees with the almost
 character indexed by the sign character $\eps$  of $W$ (see \cite[6C]{Ma95}
 and \cite{MChev}). That is 
 $$\St=|W|^{-1}\sum_{w\in W}\eps(w)R_{w}^\GG(1_{T_w}).\eqno{(\St_1)} $$
\par
(3) When $\GG=(W,\La)$, $q\in 1+\ell\ZZ_\ell$, $(|W|,\ell)=1$ and
 $\chi\in\Irr(B_0)$, this formula for $\chi(t)$ matches with that in
 Definition~\ref{def:chartab spets} by Lemma~\ref{l:conncharform} and Corollary \ref{cor:Curtis}.
\end{rem}

We conjecture that the following generalisations of
Conjectures~\ref{conj:restriction} and~\ref{c:frobspets} (for unipotent
characters) hold for these values. 

\begin{conj}   \label{conj:restriction-uni}
 Let $\GG=(W\phi,\La)$ be a simply connected $\ZZ_\ell$-spets and
 $q\in\ZZ_\ell^\times$. Let $\rho\in \Uch(\GG)$ be a unipotent character. Then
 $\rho(t)_{x=q}\in\cO$ for each $t\in S$, and regarded as an $S$-class
 function, the specialisation at $x=q$ of $\rho$ lies in $\cO[\Irr(S)]$.
\end{conj}

\begin{conj}   \label{c:frobspets-uni}
 Let $\GG=(W\phi,\La)$ be a simply connected $\ZZ_\ell$-spets and
 $q\in\ZZ_\ell^\times$. Then for each unipotent character $\rho\in\Uch(\GG)$,
 $$\Big(\sum_{t\in S/\cF}|\GG:C_{\GG}(t)|\,\rho(t)\Big)_{x=q}
       \equiv 0\pmod {|S|}$$ 
 where the sum runs over a set of $\cF$-conjugacy class representatives of elements of $S$.
\end{conj}

Conjectures \ref{conj:restriction-uni} and \ref{c:frobspets-uni} both hold when
$W$ is rational by Example \ref{exmp:DL}. We consider $\rho$ as in
Conjecture~\ref{conj:restriction-uni} as an element of $\cU(\GG)$.
If $t\in S$ is such that $\LL=C_\GG(t)$ is
$\Phi$-split for some $\Phi$ then Corollary~\ref{cor:Curtis} implies that
$\rho(t)_{x=q}\in\cO$ since $\la(1)\in\cO$ for all $\la\in\cE(\LL,X)$ and
$X\in\cM$ in the notation of that result. This is likewise true if $\ell$ does
not divide $|W_\LL|$ by Corollary~\ref{cor:val unip}. It can be checked from
the results in \cite{Ma95,BMM14} that the base change matrix between almost
characters and unipotent characters (the Fourier matrix) has denominators
only involving divisors of $|W|$. Thus, Conjecture~\ref{conj:restriction-uni}
for primes $\ell$ not dividing $|W|$ would follow from the stronger assertion:

\begin{conj}   \label{conj:almost}
 Let $\GG=(W\phi,\La)$ be a simply connected $\ZZ_\ell$-spets and
 $q\in\ZZ_\ell^\times$. Then for all $\vhi\in\Irr(W\phi)$ and all
 $\psi\in\Irr(S)$, the scalar product
 $$\bla R_\vhi,\psi\bra_S
   :=|S|^{-1}\sum_{t\in S}\psi(t)\, \big(f_\vhi^{C_\GG(t)}\big)_{x=q}$$
 is a rational integer.
\end{conj}

For the remainder of this subsection, assume that $\ell>2$ is a very good prime for $(W,\La)$ and that $(|W|,\ell)=1$. In this situation, the scalar product $\bla R_\vhi,\psi\bra_S$ is rather simple to
compute.

\begin{lem}   \label{lem:almost poly}
 Assume that $(\ell,|W|)=1$. There exist integers $b_i^{W_0}\in \ZZ$, where
 $W_0$ runs over the parabolic subgroups of $W$, such that 
 $$\bla R_\vhi,\psi\bra_T
   :=|T|^{-1}\sum_{t\in T}\psi(t)\, \left(f_\vhi^{C_\GG(t)}\right)_{x=q}
   =|T|^{-1}\sum_{W_0\le W} \psi(t) (f_\vhi^{W_0})_{x=q}\prod_{i=1}^n(\ell^a-b_i^{W_0})$$
 for all $q\in1+\ell^a\ZZ_\ell^\times$ (so $S=T \cong(\ZZ/\ell^a\ZZ)^n$),
 $\vhi\in\Irr(W)$ and $\psi\in\Irr(T)$.
\end{lem}

\begin{proof}
We reorder the sum over all $t\in T$ according to the stabiliser $W_0:=C_W(t)$
of $t$, a parabolic subgroup of $W$, to get
$$\bla R_\vhi,\psi\bra_T=|T|^{-1}\sum_{W_0\le W}m(W_0)f_\vhi^{W_0}$$
with $m(W_0)=|\{t\in T\mid C_W(t)=W_0\}|$. Now by a result of Orlik and Solomon
(see \cite[Lemma 4.7]{KMS2}), there are $b_i^{W_0}\in\ZZ$ such that
$m(W_0)=\prod_{i=1}^n(\ell^a-b_i^{W_0})$.
\end{proof}

In particular if $q=\ell^a+1$, in $\bla R_\vhi,\psi\bra_T$ we find an integral
polynomial in $\ell^a$, divided by $|T|=\ell^{an}$. In fact, we expect more.
Let $W$ be a $\ZZ_\ell$-reflection group on $\La=\ZZ_\ell^n$, with $|W|$ prime
to~$\ell$. Then for any choice of $a\ge1$ the associated Sylow $\ell$-subgroup
is $T_a=\La/\ell^a \La\cong(\ZZ/\ell^a\ZZ)^n$. Now any character
$\psi\in\Irr(T_a)$
can be considered as character of $T_b$ for all $b\ge a$ by inflation. The
following extends Conjecture~\ref{conj:almost} in this case:

\begin{conj}   \label{conj:almost 2}
 Let $\GG=(W,\La)$ be a $\ZZ_\ell$-reflection group with $(|W|,\ell)=1$.
 Then for all $a\ge1$, all $\vhi\in\Irr(W)$ and all $\psi\in\Irr(T_a)$, there
 exists a polynomial $F_{a,\vhi,\psi}\in\ZZ[y]$ with non-negative coefficients
 such that
 $$|T_b|^{-1}\sum_{t\in T_b}\psi(t)\,\big(f_\vhi^{C_\GG(t)}\big)_{x=\ell^b+1}
  = F_{a,\vhi,\psi}(\ell^b)\qquad\text{for all $b\ge a$}.$$
\end{conj}

In analogy with the spetsial case, we may denote the left hand side by
$\bla R_\vhi,\psi\bra_{T_b}$.
As a quick sanity check, we observe that $\bla R_\vhi,\psi\bra_{T_b}$
is at least a rational number: since fake degrees are rational polynomials
in~$x$, and algebraically conjugate elements of $T_b$ have the same centraliser
in $W$, any Galois automorphism permutes the summands and so the scalar product
is stable under any Galois automorphism. Furthermore if $\psi_\reg$ denotes the
regular character of $T_b$, we have
$\langle R_\vhi,\psi_\reg\bra_{T_b}=(f_\vhi^\GG)_{x=\ell^b+1}$; now
$f_\vhi^\GG$ by definition is a non-negative integer linear combination of
powers of $x$, so clearly its value at $x= \ell^b+1$ is as conjectured.

The following example illustrates Conjecture \ref{conj:almost} for the odd
dihedral groups using Lemma \ref{lem:almost poly}.

\begin{exmp}
 Let $W=G(e,e,2)$, $e\ge3$ odd, the dihedral group of order $2e$ and assume, for
 simplicity, that $q-1=\ell^a$. We have three types of elements $t\in T$, with
 centralisers $C_W(t)=W$, $W(A_1)$ and $1$, and there are $1$, $e(q-2)$,
 $(q-2)(q-e)$ such elements, respectively. The fake
 degrees $f_\varphi^\GG$ for $\varphi \in \Irr(W)$ are $1, x^e$ and
 $x^i+x^{e-i}$ for $1\le i<e/2$ and the associated (specialised) scalar
 products $\langle R_\varphi,1_T \rangle_T$ are $(q-1)^{-2}$ times
 $$(q-1)^2,\ q^e+(q-2)(e(q-1)+q),\ q^i+q^{e-i}+(q-2)(e(q-1)+2q),$$
 respectively, all divisible by $(q-1)^2$.
 For $\psi\in\Irr(T)$ non-trivial but trivial on all elements with centraliser
 $W(A_1)$ the scalar products with $\psi$ are $(q-1)^{-2}$ times
 $$0,\ q^e+eq^2-3eq+2e-1,\ q^i+q^{e-i}+eq^2-3eq+2e-2,$$
 and for $\psi\in\Irr(T)$ non-trivial on an element with centraliser $W(A_1)$
 we obtain $(q-1)^{-2}$ times
 $$0,\ q^e+(e-1)q^2-(3e-2)q+2e-2,\ q^i+q^{e-i}+(e-1)q^2-(3e-2)q+2e-3.$$
 Thus by differentiating once or twice we see that for all almost characters
 and all $\psi\in\Irr(T)$, the scalar product is integral. 
\end{exmp}

The next two examples provide further evidence for Conjecture~\ref{conj:almost}
when $\GG$ is respectively of type $G_4$ (with $\ell\equiv1\pmod 3$) and
$G_{24}$ (with $\ell\equiv 1,2,4\pmod 7$).

\begin{exmp}
Suppose $\GG$ is of type $W=G_4$, $\ell\equiv 1\pmod 3$, $q-1=\ell^a$,
$\psi=1_\TT$ and $\vhi\in\Irr(W)$ has degree $3$. Then for $t\in T$, with
$\LL(t):=C_\GG(t)$ we have
$$\begin{array}{cccc}
\LL(t) & \mbox{ \# orbits} & \vhi|_{\LL(t)} & f_\vhi^{\LL(t)} \\ \hline 1 & (q-2)(q-4)/24 & 1_\LL+1_\LL+1_\LL & 3 \\ G(1,1,3) & (q-2)/2 & 1_\LL+\omega_\LL+\omega^2|_\LL & \Phi_3 \\ W & 1 & \vhi & q^2\Phi_3\Phi_6
\end{array}$$
and we obtain
$$\langle R_\vhi,1_T\rangle_T=3(q-2)(q-4)+4(q-2)\Phi_3
  +q^2\Phi_3\Phi_6 \equiv 0 \pmod {(q-1)^2}.$$
\end{exmp}

\begin{exmp}
Suppose $\GG$ is of type $W=G_{24}$, $\ell\equiv 1,2,4\pmod 7$, $q-1=\ell^a$,
$\psi=1_\TT$ and $\vhi=\vhi_{3,1}\in\Irr(W)$ has degree $3$. Then we have:
$$\begin{array}{cccc}
\LL(t) & \mbox{ \# orbits} & \vhi|_{\LL(t)} & f_\vhi^{\LL(t)} \\
 \hline
 1 & (q-2)(q-10)(q-12)/336 & 3 \cdot 1_\LL & 3 \\
 G(2,1,1) & (q-2)(q-8)/8 & 2 \cdot 1_\LL+ \mbox{sgn}_\LL & q+2 \\
 G(1,1,3) & (q-2)/2 & 1_\LL+\mbox{perm}_\LL & 1+q\Phi_2 \\
 G(2,1,2) & (q-2)/2 & 1_\LL+\mbox{sq}_\LL & 1+q\Phi_4 \\
 W & 1 & \phi_{3,1} & (1-q^2+q^6)q\Phi_3\Phi_6
\end{array}$$
and we obtain
$$\langle R_\vhi,1_T\rangle_T=(q-1)^3 (q^8 + 3q^7 + 7q^6 + 13q^5 + 21q^4
  + 31q^3 + 43q^2 + 78q + 146). $$
\end{exmp}

\subsection{The $\ZZ_2$-spets $G_{24}(q)$ related to the Benson--Solomon fusion system $\Sol(q)$}   \label{s:sol}

Here we give evidence for
Conjectures~\ref{conj:restriction-uni} and \ref{c:frobspets-uni} by considering
the $\ZZ_2$-spets $G_{24}(q)$ when $q\equiv 1\pmod 4$ (the case of
$q\equiv 3\pmod 4$ is Ennola dual and entirely similar so we do not include it
here). Note that $2$ divides the order of $W$ and in fact, $2$ is a bad prime
for $W$ in the sense of \cite[\S 2.4]{KMS}.

Set $\GG=(W,\La)$ where $\La=\ZZ_2^3$, and $W\le\GL(\La)$ is the
$\ZZ_2$-reflection group $G_{24}$. Assume $q\in1+4\ZZ_2$. Let $\GG(q)$ be
the associated $\ZZ_2$-spets, to which we may associate a $2$-fusion
system~$\cF$ exactly as in \cite[\S 3.1]{KMS}, by taking $q$-homotopy fixed
points of the associated $2$-compact group BDI$(4)$. It is a result of Benson
(see \cite[Thm~4.5]{LO02}) that $\cF=\Sol(q)$ is the Benson--Solomon $2$-fusion
system on a Sylow $2$-subgroup~$S$ of $\Spin_7(q)$. Note that $|S| = 1024\,l^3$
where $l:=2^{f-3}$ with $f:=\nu_2(q^2-1)$. It is well known that $\cF$ does not
occur as the fusion system of any block of any finite group (see
\cite[Thm~1.1]{K06} for the case $q=3$ and \cite[Prop.~9.34]{C11} for
general~$q$). 

The third author in \cite{S19} used an adhoc method to define the `type' of
$C_\GG(s)$ for each fully $\cF$-centralised element $s$ of $S$, where he also
determined the number of $\cF$-classes of each type (see \cite[Tab.~10]{S19}).
We expect that this method fits within the general framework of \cite{KMS},
but we have not been able to prove this.

Using the formula \eqref{e:uniformula}, the values $\chi(s)$ taken by unipotent
characters $\chi\in\Uch(\GG)$ at $s\in S$ may now be obtained using
{\sf Chevie} \cite{MChev}, and we list these in Table~\ref{t:G24values}.

\begin{table}[htb]
\caption{Unipotent character table of $G_{24}(q)$, $q\equiv1\pmod4$}
\label{t:G24values}
$$\begin{array}{c|ccc}
 \# \mbox{ $\cF$-cl.} & 1 & 1 & 1\\
 C_\GG(s)(q)& \GG(q)& B_3(q)& A_3(q)\\
\hline
 \phi_{1,0}&                         1& 1& 1\\
 \phi_{3,1}&  \frac{\sq7}{14}q\Ph3\Ph4\Ph6\Ph7'\Ph{14}''& \hlf\Ph3\Ph4\Ph6& \hlf\Ph3\Ph4\\
 \phi_{3,3}& -\frac{\sq7}{14}q\Ph3\Ph4\Ph6\Ph7''\Ph{14}'& \hlf\Ph3\Ph4\Ph6& \hlf\Ph3\Ph4\\
 \phi_{6,2}&  \hlf q\Ph2^2\Ph3\Ph6\Ph{14}& \hlf\Ph2^2\Ph3\Ph6& \hlf\Ph2^2\Ph3\\
    B_2\co2&  \hlf q\Ph1^2\Ph3\Ph6\Ph7& -\hlf\Ph1^2\Ph3\Ph6& -\hlf\Ph1^2\Ph3\\
 \phi_{7,3}&  q^3\Ph7\Ph{14}& q(q^6\pl2q^4\pl3q^2\pl1)& (q^5\pl q^4\pl2q^3\pl q^2\pl q\pl1)\\
 \phi_{8,4/5}&  \hlf q^4\Ph2^3\Ph4\Ph6\Ph{14}& \hlf q\Ph2^3\Ph4\Ph6& q\Ph2^2\Ph4\\
  G_{24}[\pm I]&  \hlf q^4\Ph1^3\Ph3\Ph4\Ph7& -\hlf q\Ph1^3\Ph4\Ph3& .\\
 \phi_{7,6}&  q^6\Ph7\Ph{14}& q^2(q^6\pl3q^4\pl2q^2\pl1)& q(q^5\pl q^4\pl q^3\pl2q^2\pl q\pl1)\\
 \phi_{3,8}&-\frac{\sq7}{14}q^8\Ph3\Ph4\Ph6\Ph7''\Ph{14}'& \hlf q^3\Ph3\Ph4\Ph6& \hlf q^2\Ph3\Ph4\\
\phi_{3,10}& \frac{\sq7}{14}q^8\Ph3\Ph4\Ph6\Ph7'\Ph{14}''& \hlf q^3\Ph3\Ph4\Ph6& \hlf q^2\Ph3\Ph4\\
 \phi_{6,9}&  \hlf q^8\Ph2^2\Ph3\Ph6\Ph{14}& \hlf q^3\Ph2^2\Ph3\Ph6& \hlf q^2\Ph2^2\Ph3\\
  B_2\co1^2&  \hlf q^8\Ph1^2\Ph3\Ph6\Ph7& -\hlf q^3\Ph1^2\Ph3\Ph6& -\hlf q^2\Ph1^2\Ph3\\
\phi_{1,21}&  q^{21}& q^9& q^6\\
\end{array}$$
\end{table}

\begin{table}[htb]
$$\begin{array}{c|cccccccc}
  & 2l-1& 2(l-1)& l-1& l\\
    & B_2(q)\Ph1& A_2(q)\Ph1& A_1(q)^2\Ph1& A_1(q^2)\Ph1\\
\hline
   \phi_{1,0}&  1& 1& 1& 1\\
 \phi_{3,1/3}&  \hlf(q+2)\Ph4& \Ph3& \hlf(q^2\pl2q\pl3)& \hlf\Ph4\\
   \phi_{6,2}&  \hlf(3q^2\mn q\pl4)\Ph2& 2\Ph3& \hlf(q\pl5)\Ph2& \hlf(q^2\pl3)\\
      B_2\co2&  \hlf q\Ph1^2& .& -\hlf\Ph1^2& -\hlf\Ph1\Ph2\\
   \phi_{7,3}&  2\Ph2\Ph4-1& 2\Ph2\Ph4\mn q^3& q^2\pl4q\pl2& q^2\\
 \phi_{8,4/5}&  \Ph2^2\Ph4& \Ph2^3& 2\Ph2^2& .\\
   \phi_{7,6}&  2\Ph2\Ph4-q^3& 2\Ph2\Ph4\mn1& 2q^2\pl4q\pl1& -1\\
 \phi_{3,8/10}&  \hlf q(2q+1)\Ph4& q\Ph3& \hlf(3q^2\pl2q\pl1)& -\hlf\Ph4\\
   \phi_{6,9}&  \hlf q(4q^2\mn q\pl3)\Ph2& 2q\Ph3& \hlf\Ph2(5q\pl1)& -\hlf(3q^2\pl1)\\
    B_2\co1^2&  \hlf q\Ph1^2& .& -\hlf\Ph1^2& -\hlf\Ph1\Ph2\\
  \phi_{1,21}&  q^4& q^3& q^2& -q^2\\
\end{array}$$
\end{table}

\begin{table}[htb]
$$\begin{array}{c|cccccccc}
  &  (l\mn1)(2l\mn3)& l& \frac{1}{21}(4l\mn9)(l\mn1)(l\mn2)& l(l\mn1)\\
  & A_1(q)\Ph1^2& A_1(q)\Ph1\Ph2& \Ph1^3& \Ph1^2\Ph2\\
\hline
 \phi_{1,0}& 1& 1& 1& 1\\
 \phi_{3,1/3}& q+2& 1& 3& 1\\
 \phi_{6,2}& 2(q\pl2)& \Ph2& 6& 2\\
  B_2\co2& .& \Ph1& .& .\\
 \phi_{7,3}& 3q+4& q& 7& 1\\
 \phi_{8,4/5}& 4\Ph2& .& 8& .\\
 \phi_{7,6}& 4q+3& -1& 7& -1\\
 \phi_{3,8/10}& 2q+1& -q& 3& -1\\
 \phi_{6,9}& 2(2q\pl1)& -\Ph2& 6& -2\\
 B_2\co1^2&. & \Ph1& .& .\\
\phi_{1,21}& q& -q& 1& -1\\
\end{array}$$
\end{table}

Note that the other six (cuspidal) unipotent characters are of $2$-defect $0$
and vanish on all $2$-elements $s\ne1$, as do $G_{24}[\pm I]$ on all classes
in the second and third part of the table. It is now straightforward to check
that:

\begin{prop}   \label{prop:solq1}
 Conjecture \ref{c:frobspets-uni} holds for all unipotent
 characters $\chi\in\Uch(\GG)$.
\end{prop} 

Moreover, when $l=1$, the \texttt{MAGMA} \cite{BCP97} package \cite{PS21}
developed jointly by the third author and Parker may be used to provide
explicit representatives for each of the six $\cF$-classes of~$S$. Using this
information, we observe that $\langle \chi,\psi\rangle\in\cO$ for each
$\psi \in \Irr(S)$, and conclude:

\begin{prop}   \label{prop:solq2}
 Conjecture \ref{conj:restriction-uni} holds for
 $\chi\in\Uch(\GG)$ whenever $\nu_2(q^2-1)=3$.
\end{prop}

Naturally, we expect the hypothesis on $q$ in Proposition \ref{prop:solq2} to
be unnecessary but we did not check this.

\subsection{The $\ZZ_5$-spets $G_{29}(q)$}   \label{s:agu}
We have also computed the character values for the $\ZZ_5$-spets $G_{29}(q)$
related to the Aguad\'{e} $5$-compact group $\bX_{29}$ when $q\equiv 1\pmod 5$.
Thus $5$ divides the order of $W$ and is a bad prime for $W$ in the sense of
\cite[\S 2.4]{KMS}. 

Set $\GG=(W,\La)$ where $\La=\ZZ_5^4$, and $W\le\GL(\La)$ is the
$\ZZ_5$-reflection group $G_{29}$.  Suppose that $q \in 1+ 5^f\ZZ_5^\times$.
Then by \cite[Thm~3.2]{KMS} there is an exotic fusion system $\cF=\cG_{29}$
associated to the pair $(\GG,q)$ on a group $S$ isomorphic to a Sylow
$5$-subgroup of $\SL_5(q)$ of order $5^{4f+1}$ obtained by taking $q$-homotopy
fixed points of the associated Aguad\'{e} $5$-compact group $\bX_{29}$. In fact
$\cG_{29}$ occurs naturally as a subsystem of the $5$-fusion system of
$E_8(2^{5^{f-1}})$ but we shall not need this. Using \cite[Lemma 8.2]{KMS}, we
see that every element of $S$ is $\cF$-conjugate to an element of the maximal
$5$-torus $T$, and so we may choose a set $S^\cF \subseteq T$ of fully
$\cF$-centralised class representatives. Moreover by the discussion in
\cite[\S 6.2]{KMS} there is a well-defined centraliser spets
$C_{\GG(q)}(s)=(W(s),\Lambda)$ for each $s \in S^\cF$.  Using the formula
\eqref{e:uniformula}, the values $\chi(s)$ taken by unipotent characters
$\chi\in\Uch(\GG)$ at $s\in S$ may now be obtained using
{\sf Chevie} \cite{MChev}, and these are listed in Table~\ref{t:G29values}. 

In \cite{CV25}, the authors prove that the total number of $\cF$-classes of
elements in $S$ is a polynomial
\begin{equation}   \label{eq:cv}
  \frac{1}{7680} \left( l^4 + 40 l^3 + 530 l^2 + 2720 l + 5925 \right)
\end{equation}
where $l=5^f$. In fact, via an elementary inclusion/exclusion argument, the
number of $\cF$-classes of elements in $S$ with a given $W$-stabiliser is
necessarily polynomial in $l$. We have used this fact, together with
(\ref{eq:cv}) to interpolate these numbers. Moreover, when $l=1$, it is
straightforward to obtain explicit representatives for each of the five
$\cF$-classes of~$S$. We observe that $\langle \chi,\psi\rangle\in\cO$ for
each $\psi\in \Irr(S)$.

\begin{prop}   \label{prop:g29}
 Let $\GG=G_{29}(q)$. Then we have:
 \begin{enumerate}[\rm(a)]
  \item Conjecture \ref{c:frobspets-uni} holds for all unipotent characters
   $\chi\in\Uch(\GG)$;
  \item  Conjecture \ref{conj:restriction-uni} holds for $\chi\in\Uch(\GG)$
   whenever $\nu_5(q-1)=1$.
 \end{enumerate}
\end{prop}

As for Proposition \ref{prop:solq2}, we expect the hypothesis on $q$ in
Proposition \ref{prop:g29}(b) to be unnecessary but we did not check this.

\tiny
\begin{table}[htb]
\caption{Unipotent character table of $G_{29}(q)$, $q \equiv 1 \pmod 5$}\label{t:G29values}

\begin{longtable}{c|ccc}

 \# $\cF$\mn cl & 1 & 1& \\
 $C_{\GG}(s)(q)$ & $\GG(q)$ & $A_4(q)$ & \\
\hline

$\phi_{1,0}$ &
1 & 1& \\ 

$\phi_{4,1} $ & $
\frac{1\pl i}{4} q \Ph2^2 \Ph6 \Ph{10} \Ph{20}'' \Ph{12}'' \Ph8'' (\Ph4')^2$ & $
q \Ph2 \Ph4 $ & \\  

$\phi_{4,3} $ & $
\frac{1\mn i}{4} q \Ph2^2 \Ph6 \Ph{10} \Ph{20}' \Ph{12}' \Ph8' (\Ph4'')^2 $ & $
q \Ph2 \Ph4 $ & \\  

$\phi_{4,4} $ & $
\frac12 q \Ph4^3 \Ph{12}\Ph{20} $ & $
q \Ph2\Ph4 $ & \\  

$\phi_{10,2} $ & $
q^2 \Ph5\Ph8\Ph{10}\Ph{20} $ & $
\Ph4\Ph5 $ & \\  

$\phi_{16,3} $ & $
\frac12 q^3 \Ph4^4 \Ph8\Ph{12}\Ph{20} $ & $
\Ph2\Ph4(q^4\pl q^3\pl q^2\pl 1) $ &\\  

$\phi_{16,5} $ & $
\frac12 q^3 \Ph4^4 \Ph8\Ph{12}\Ph{20} $ & $
\Ph2\Ph4(q^4\pl q^3\pl q^2\pl 1) $ &\\   

$\phi_{15,4}'' $ & $
\frac12 q^4 \Ph3\Ph5\Ph6\Ph8\Ph{10}\Ph{20} $ & $
q\Ph3\Ph5 $ & \\

$\phi_{5,8} $ & $
\frac12 q^4 \Ph5\Ph8\Ph{10}\Ph{12}\Ph{20} $ & $
q^2\Ph5 $ & \\

$\phi_{10,6} $ & $
\frac12 q^4 \Ph4^2\Ph5\Ph{10}\Ph{12}\Ph{20} $ & $
q\Ph4\Ph5 $ & \\

$\phi_{15,4}' $ & $
q^4 \Ph3\Ph5\Ph6\Ph{10}\Ph{12}\Ph{20} $ & $
\Ph3\Ph5\Ph6 $ &\\

$\phi_{20,5} $ & $
\frac{1\mn i}{4} q^5 \Ph2^2\Ph5\Ph6\Ph{10}\Ph{20}\Ph{12}''\Ph8'(\Ph4')^2 $ & $
q\Ph2\Ph4\Ph5 $ &\\

$\phi_{20,7} $ & $
\frac{1\pl i}{4} q^5 \Ph2^2\Ph5\Ph6\Ph{10}\Ph{20}\Ph{12}'\Ph8''(\Ph4'')^2 $ & $
q\Ph2\Ph4\Ph5 $ &\\

$\phi_{20,6} $ & $
\frac12 q^5 \Ph4^3\Ph5\Ph{10}\Ph{12}\Ph{20} $ & $
q\Ph2\Ph4\Ph5 $ &\\

$\phi_{24,6} $ & $
\frac{1}{20} q^6 \Ph2^4\Ph3\Ph5\Ph6\Ph8\Ph{12}\Ph{20} $ & $
\frac15 \Ph3\Ph5\Ph4\Ph2^2 $ &\\

$\phi_{6,10}' $ & $
\frac15 q^6 \Ph3\Ph5\Ph6\Ph8\Ph{10}\Ph{12}\Ph{20} $ & $
\frac15 \Ph3\Ph5(4q^4\mn 7q^3\pl 3q^2\pl 3q\mn 1) $ &\\

$\phi_{6,10}'' $ & $
\frac15 q^6 \Ph3\Ph5\Ph6\Ph8\Ph{10}\Ph{12}\Ph{20} $ & $
\mn \frac15 \Ph3\Ph5(q^4\mn 3q^3\mn 3q^2\pl 7q\mn 4) $ &\\

$\phi_{6,10}''' $ & $
\mn \frac{1}{20} q^6 \Ph3\Ph5\Ph6\Ph8\Ph{10}\Ph{12}\Ph{20}''(\Ph4'')^4 $ & $
\mn \frac15 \Ph3\Ph4\Ph5(q^2\mn 3q\pl 1) $ & \\

$\phi_{6,10}'''' $ & $
\mn \frac{1}{20} q^6 \Ph3\Ph5\Ph6\Ph8\Ph{10}\Ph{12}\Ph{20}'(\Ph4')^4 $ & $
\mn \frac15 \Ph3\Ph4\Ph5(q^2\mn 3q\pl 1) $ &\\

$\phi_{24,9} $ & $
\frac14 q^6 \Ph2^2\Ph3\Ph4^3\Ph6\Ph{10}\Ph{12}\Ph{20} $ & $
q\Ph2^2\Ph3\Ph4\Ph6 $ & \\

$\phi_{24,7} $ & $
\frac14 q^6 \Ph2^2\Ph3\Ph4^3\Ph6\Ph{10}\Ph{12}\Ph{20} $ & $
q\Ph2^2\Ph3\Ph4\Ph6 $ &\\

$\phi_{30,8} $ & $
\frac14 q^6 \Ph3\Ph4^2\Ph5\Ph6\Ph8\Ph{12}\Ph{20} $ & $
q\Ph3\Ph5\Ph4 $ & \\

$\phi_{6,12} $ & $
\frac14 q^6 \Ph3\Ph4^2\Ph6\Ph8\Ph{10}\Ph{12}\Ph{20} $ & $
q^3\Ph3\Ph4 $ &\\

$G_{29}[1] $ & $
\frac{1}{20} q^6 \phi_1^4\Ph3\Ph6\Ph8\Ph{10}\Ph{12}\Ph{20} $ & $
\frac15 \phi_1^4\Ph2^2\Ph3\Ph4 $ &\\

$G_{29}[\ze_5] $ & $
\frac15 q^6 \phi_1^4\Ph2^4\Ph3\Ph4^4\Ph6\Ph8\Ph{12} $ & $
\mn \frac15 \phi_1^4\Ph2^2\Ph3\Ph4 $ & \\

$G_{29}[\ze_5^2] $ & $
\frac15 q^6 \phi_1^4\Ph2^4\Ph3\Ph4^4\Ph6\Ph8\Ph{12} $ & $
\mn \frac15 \phi_1^4\Ph2^2\Ph3\Ph4 $ & \\

$G_{29}[\ze_5^3] $ & $
\frac15 q^6 \phi_1^4\Ph2^4\Ph3\Ph4^4\Ph6\Ph8\Ph{12} $ & $
\mn \frac15 \phi_1^4\Ph2^2\Ph3\Ph4 $ &\\

$G_{29}[\ze_5^4] $ & $
\frac15 q^6 \phi_1^4\Ph2^4\Ph3\Ph4^4\Ph6\Ph8\Ph{12} $ & $
\mn \frac15 \phi_1^4\Ph2^2\Ph3\Ph4 $ &\\

$\phi_{20,9} $ & $
\frac{1\pl i}{4} q^9 \Ph2^2\Ph5\Ph6\Ph{10}\Ph{20}\Ph{12}'\Ph8''(\Ph4'')^2 $ & $
q^2\Ph2\Ph4\Ph5 $ &\\

$\phi_{20,11} $ & $
\frac{1\mn i}{4} q^9 \Ph2^2\Ph5\Ph6\Ph{10}\Ph{20}\Ph{12}''\Ph8'(\Ph4')^2 $ & $
q^2\Ph2\Ph4\Ph5 $ &\\

$\phi_{20,10} $ & $
\frac12 q^9 \Ph4^3\Ph5\Ph{10}\Ph{12}\Ph{20} $ & $
q^2\Ph2\Ph4\Ph5 $ &\\

$\phi_{15,12}' $ & $
q^{12}\Ph3\Ph5\Ph6\Ph{10}\Ph{12}\Ph{20} $ & $
q^2\Ph3\Ph5\Ph6 $ & \\

$\phi_{15,12}'' $ & $
\frac12 q^{12}\Ph3\Ph5\Ph6\Ph8\Ph{10}\Ph{20} $ & $
q^3\Ph3\Ph5 $ & \\

$\phi_{5,16} $ & $
\frac12 q^{12}\Ph5\Ph8\Ph{10}\Ph{12}\Ph{20} $ & $
q^4\Ph5 $ & \\

$\phi_{10,14} $ & $
\frac12 q^{12}\Ph4^2\Ph5\Ph{10}\Ph{12}\Ph{20} $ & $
q^3\Ph4\Ph5 $ & \\

$\phi_{16,13} $ & $
\frac12 q^{13}\Ph4^4\Ph8\Ph{12}\Ph{20} $ & $
q^3\Ph2\Ph4(q^4\pl q^2\pl q\pl 1) $ & \\

$\phi_{16,15} $ & $
\frac12 q^{13}\Ph4^4\Ph8\Ph{12}\Ph{20} $ & $
q^3\Ph2\Ph4(q^4\pl q^2\pl q\pl 1) $ & \\

$\phi_{10,18} $ & $
q^{18}\Ph5\Ph8\Ph{10}\Ph{20} $ & $
q^4\Ph4\Ph5 $ & \\

$\phi_{4,21} $ & $
\frac{1\mn i}{4} q^{21}\Ph2^2\Ph6\Ph{10}\Ph{20}'\Ph{12}'\Ph8'(\Ph4'')^2 $ & $
q^6\Ph2\Ph4 $ & \\

$\phi_{4,23} $ & $
\frac{1\pl i}{4} q^{21}\Ph2^2\Ph6\Ph{10}\Ph{20}''\Ph{12}''\Ph8''(\Ph4')^2 $ & $
q^6\Ph2\Ph4 $ & \\

$\phi_{4,24} $ & $
\frac12 q^{21}\Ph4^3\Ph{12}\Ph{20} $ & $
q^6\Ph2\Ph4 $ & \\

$\phi_{1,40} $ & $
q^{40} $ & $ q^{10} $ & \\   

\end{longtable}
\end{table}
\newpage

\begin{longtable}{c|cc}
 
 \# $\cF$\mn cl & $\frac{1}{4}(l\mn 1)$& \\
 $C_\GG(s)(q)$ & $\GG_{4,4,3}(q)\Phi_1$ & \\  \hline

$\phi_{1,0}$ &  1\\  

$\phi_{4,1} $& $
\mn \frac{1\pl i}{4}(2i q^4\mn q^5\mn q^4\mn 3q^2\mn 2\mn 2i\mn q)\Ph8''$ \\  

$\phi_{4,3} $ &  $
\frac{1\mn i}{4}(2i q^4\pl q^5\pl q^4\pl 3q^2\pl 2\mn 2i\pl q)\Ph8' $ \\  

$\phi_{4,4}  $& $\frac12 \Ph2\Ph4 (q^4\mn q\pl 2) $ \\  

$\phi_{10,2} $ & $(q^3\pl 2q^2\pl q\pl 1)\Ph8 $\\  

$\phi_{16,3} $& $\Ph2\Ph4^2\Ph8 $\\  

$\phi_{16,5} $ &  $\Ph2\Ph4^2\Ph8 $\\

$\phi_{15,4}'' $ & $\frac12 \Ph3\Ph8(2q^3\pl 2q^2\mn q\pl 2) $ \\

$\phi_{5,8} $ & $\frac12 q\Ph8(2q^3\pl q^2\pl q\pl 1) $ \\

$\phi_{10,6} $ &  $
q^9\pl q^8\pl \frac32 q^7\pl \frac12 q^6\pl \frac32 q^5\pl q^4\pl \frac32 q^3\pl \frac12 q^2\pl \frac12 q\pl 1 $ \\

$\phi_{15,4}' $ & $ q\Ph3(q^6\pl q^4\pl q^3\pl q\pl 1) $\\

$\phi_{20,5} $ & $
\frac{1\mn i}{4}q\Ph8'(\mn 2iq^2\pl q^8\pl 2iq^4\pl q^7\mn 2iq^3\pl 2q^6\pl 2iq^7\pl 7q^5\mn 4iq\pl 7q^4\pl 2iq^6\pl 10q^3\pl 4iq^5\pl 6q^2\pl 2\mn 2i\pl 4q) $ \\

$\phi_{20,7} $ &  $
\mn \frac{1\pl i}{4}q\Ph8''(\mn 2iq^2\mn q^8\pl 2iq^4\mn q^7\mn 2iq^3\mn 2q^6\pl 2iq^7\mn 7q^5\mn 4iq\mn 7q^4\pl 2iq^6\mn 10q^3\pl 4iq^5\mn 6q^2\mn 2\mn 2i\mn 4q) $\\

$\phi_{20,6} $ & $\frac12 q\Ph2\Ph4(q^7\pl 2q^5\pl q^4\pl 2q^3\pl 2q\pl 2) $\\

$\phi_{24,6} $ & $\frac12 q\Ph3\Ph4\Ph8\Ph2^2 $\\

$\phi_{6,10}' $& $q^3\Ph3\Ph8 $\\

$\phi_{6,10}'' $ &  $q^3\Ph3\Ph8 $\\

$\phi_{6,10}''' $ & $ \frac{1\pl i}{4} q\Ph3\Ph8''(\Ph8')^2(\Ph4'')^2 $\\

$\phi_{6,10}'''' $ & $\frac{1\mn i}{4} q\Ph3\Ph8'(\Ph8'')^2(\Ph4')^2 $\\

$\phi_{24,9} $ & $
\mn \left(\frac{3}{20}\pl \frac{i}{20}\right)\Ph3\Ph4 q
(4iq^5\mn 5q^6\pl 3iq^4\mn 2q^5\mn 6iq^3\mn 4q^4\pl 3iq^2\mn 2q^3\pl 4iq\mn 4q^2\mn 2q\mn 5) $\\

$\phi_{24,7} $ & $
\mn \left(\mn \frac{3}{20}\pl \frac{i}{20}\right)\Ph3\Ph4 q
(4iq^5\pl 5q^6\pl 3iq^4\pl 2q^5\mn 6iq^3\pl 4q^4\pl 3iq^2\pl 2q^3\pl 4iq\pl 4q^2\pl 2q\pl 5) $\\

$\phi_{30,8} $ & $ q\Ph5\Ph3\Ph8 $\\

$\phi_{6,12} $ & $ \frac12 q\Ph3\Ph8^2 $\\

$G_{29}[1] $ & $ 0 $\\

$G_{29}[\ze_5] $ & $ 0 $\\

$G_{29}[\ze_5^2] $ & $0 $\\

$G_{29}[\ze_5^3] $ & $0 $\\

$G_{29}[\ze_5^4] $ & $0 $\\

$\phi_{20,9} $ & $
\frac14 q\Ph8''(4q^8\pl 4iq^6\pl 8q^7\pl 8iq^5\pl 8q^6\pl 9iq^4\pl 12q^5\pl 11iq^3\pl 5q^4\pl 4iq^2\pl 3q^3\pl 3iq\pl 1\pl i\mn q) $\\

$\phi_{20,11} $ &  $
\mn \frac14 q\Ph8'( \mn 4q^8\pl 4iq^6\mn 8q^7\pl 8iq^5\mn 8q^6\pl 9iq^4\mn 12q^5\pl 11iq^3\mn 5q^4\pl 4iq^2\mn 3q^3\pl 3iq\mn 1\pl i\pl q) $\\

$\phi_{20,10} $ &  $
\frac12 q\Ph2\Ph4(2q^7\pl 2q^6\pl 2q^4\pl q^3\pl 2q^2\pl 1) $\\

$\phi_{15,12}' $ &  $ q^3\Ph3(q^6\pl q^5\pl q^3\pl q^2\pl 1) $\\

$\phi_{15,12}'' $ & $ \frac12 q^3\Ph3\Ph8(2q^3\mn q^2\pl 2q\pl 2) $\\

$\phi_{5,16} $ & $\frac12 q^4\Ph8(q^3\pl q^2\pl q\pl 2) $\\

$\phi_{10,14} $ & $
\frac12 q^3(2q^9\pl q^8\pl q^7\pl 3q^6\pl 2q^5\pl 3q^4\pl q^3\pl 3q^2\pl 2q\pl 2) $\\

$\phi_{16,13} $ &  $q^3\Ph2\Ph4^2\Ph8 $ \\

$\phi_{16,15} $ & $q^3\Ph2\Ph4^2\Ph8 $\\

$\phi_{10,18} $ & $ q^5\Ph8(q^3\pl q^2\pl 2q\pl 1) $\\

$\phi_{4,21} $ & $
\mn \frac14 q^5(iq^4\mn 4q^5\pl 3iq^3\mn q^4\mn 3q^3\pl 3iq\mn 1\pl i\pl q)\Ph8' $\\

$\phi_{4,23} $ & $
\frac14 q^5(iq^4\pl 4q^5\pl 3iq^3\pl q^4\pl 3q^3\pl 3iq\pl 1\pl i\mn q)\Ph8'' $ \\

$\phi_{4,24} $ & $ \frac12 q^5\Ph2\Ph4(2q^4\mn q^3\pl 1) $\\

$\phi_{1,40} $ & $ q^{12} $\\   

\end{longtable}
\newpage

\newpage
$\begin{array}{c|ccccc}
 
 \# \cF\mn cl & \frac{1}{4}(l\mn 1)& \frac{1}{4}(l\mn 1) & \frac{1}{4}(l\mn 5) & \frac{1}{4}(l\mn 5) &\\ 
 C_\GG(s)(q) & B_3(q) \Phi_1 & A_3(q)\Phi_1 & A_3(q)\Phi_1 & A_2(q)A_1(q)\Phi_1 \\
\hline
 
\phi_{1,0} & 1 & 1 & 1 & 1 \\   

\phi_{4,1} &  \frac12 \Ph2 \Ph6 (q^2\pl q\pl 2)& \Ph2\Ph4& \Ph2\Ph4& \Ph2^2\\  

\phi_{4,3}& \frac12 \Ph2 \Ph6 (q^2\pl q\pl 2)& \Ph2\Ph4& \Ph2\Ph4& \Ph2^2\\  

\phi_{4,4}& q^5\pl q^3\pl q\pl 1& \Ph2\Ph4& \Ph2\Ph4& \Ph2^2\\  

\phi_{10,2}& q^6\pl q^5\pl 2q^4\pl q^3\pl 2q^2\pl q\pl 2& \Ph4(q^2\pl 2q\pl 2)&
 \Ph4(q^2\pl 2q\pl 2)& q^3\pl 3q^2\pl 3q\pl 3 \\  

\phi_{16,3} &  \Ph4(q^5\pl q^4\pl 2q^3\pl q^2\pl 2q\pl 1)& 
 \Ph2\Ph4(q^2\pl q\pl 2)& \Ph2\Ph4(q^2\pl q\pl 2)& \Ph2(3q^2\pl 2q\pl 3)\\  

\phi_{16,5}& \Ph4(q^5\pl q^4\pl 2q^3\pl q^2\pl 2q\pl 1)& 
 \Ph2\Ph4(q^2\pl q\pl 2)& \Ph2\Ph4(q^2\pl q\pl 2)& \Ph2(3q^2\pl 2q\pl 3)\\

\phi_{15,4}''& \frac12 \Ph3\Ph6(q^4\pl q^3\pl 3q^2\pl 3q\pl 2)& 
 \Ph3(q^3\pl q^2\pl 2q\pl 1)& \Ph3(q^3\pl q^2\pl 2q\pl 1)& (3q\pl 2)\Ph3\\

\phi_{5,8}& \frac12 q^8\pl \frac12 q^7\pl q^5\pl q^4\pl \frac12 q^2\pl \frac12 q\pl 1& 
q(q^3\pl q^2\pl 2q\pl 1)& q(q^3\pl q^2\pl 2q\pl 1)& q^3\pl 2q^2\pl q\pl 1\\

\phi_{10,6}& \frac12 q(q^7\mn q^6\pl 4q^5\pl 2q^4\pl 4q^3\pl 4q^2\pl 3q\pl 3)&  
\Ph4(q^3\pl q^2\pl 2q\pl 1)& \Ph4(q^3\pl q^2\pl 2q\pl 1)& 2q^3\pl 3q^2\pl 4q\pl 1\\

\phi_{15,4}'& \Ph3\Ph6(q^3\pl 2q^2\pl q\pl 1)& \Ph3(q^3\pl 2q^2\pl 2)&
 \Ph3(q^3\pl 2q^2\pl 2)&  \Ph3(q^2\pl q\pl 3)\\

\phi_{20,5}& \frac12 \Ph2\Ph6(2q^5\pl 4q^4\pl 4q^3\pl 5q^2\pl 3q\pl 2) &
 \Ph2\Ph4(2q^2\pl 2q\pl 1)& \Ph2\Ph4(2q^2\pl 2q\pl 1)& (q^2\pl 2q\pl 2)\Ph2^2\\

\phi_{20,7}&  \frac12 \Ph2\Ph6(2q^5\pl 4q^4\pl 4q^3\pl 5q^2\pl 3q\pl 2)& 
\Ph2\Ph4(2q^2\pl 2q\pl 1)& \Ph2\Ph4(2q^2\pl 2q\pl 1)& (q^2\pl 2q\pl 2)\Ph2^2\\

\phi_{20,6}& q^8\pl 2q^7\pl 2q^6\pl 4q^5\pl 3q^4\pl 3q^3\pl 2q^2\pl 2q\pl 1& 
\Ph2\Ph4(2q^2\pl 2q\pl 1)& \Ph2\Ph4(2q^2\pl 2q\pl 1)& (q^2\pl 2q\pl 2)\Ph2^2\\

\phi_{24,6}& q\Ph3\Ph6\Ph2^3& \Ph4\Ph3\Ph2^2& \Ph4\Ph3\Ph2^2&  2\Ph3\Ph2^2\\

\phi_{6,10}' & q^2\Ph2\Ph3\Ph6&  \Ph3(q^3\pl q^2\mn q\pl 1)& 
 q\Ph3(q^3\mn q^2\pl q\pl 1)& \Ph3\Ph4\\

\phi_{6,10}''& q^2\Ph2\Ph3\Ph6&  q\Ph3(q^3\mn q^2\pl q\pl 1)& 
 \Ph3(q^3\pl q^2\mn q\pl 1)& \Ph3\Ph4\\

\phi_{6,10}'''& \frac12 q\Ph2\Ph3\Ph4\Ph6& q\Ph3\Ph4& q\Ph3\Ph4& 2q\Ph3\\

\phi_{6,10}''''&  \frac12 q\Ph2\Ph3\Ph4\Ph6& q\Ph3\Ph4& q\Ph3\Ph4& 2q\Ph3\\

\phi_{24,9}& \frac12 q\Ph2\Ph3\Ph6(3q^2\pl 2q\pl 3)& \Ph2^2\Ph3\Ph4&
 \Ph2^2\Ph3\Ph4& 2\Ph2^2\Ph3\\

\phi_{24,7}& \frac12 q\Ph2\Ph3\Ph6(3q^2\pl 2q\pl 3)& \Ph2^2\Ph3\Ph4&
 \Ph2^2\Ph3\Ph4& 2\Ph2^2\Ph3\\

\phi_{30,8}& q\Ph2\Ph3(2q^2\pl q\pl 2)\Ph6& \Ph3(q^2\pl 3q\pl 1)\Ph4&
 \Ph3(q^2\pl 3q\pl 1)\Ph4& 2\Ph3(q^2\pl 3q\pl 1)\\

\phi_{6,12}& q\Ph2\Ph3\Ph6^2& q\Ph3\Ph4& q\Ph3\Ph4& 2q\Ph3\\

G_{29}[1]&  0& 0& 0& 0 \\

G_{29}[\ze_5]& 0& 0& 0& 0\\

G_{29}[\ze_5^2]& 0& 0& 0& 0\\

G_{29}[\ze_5^3]& 0& 0& 0& 0\\

G_{29}[\ze_5^4]& 0& 0& 0& 0\\

\phi_{20,9}& \frac12 q\Ph2\Ph6(2q^5\pl 3q^4\pl 5q^3\pl 4q^2\pl 4q\pl 2)& 
 q\Ph2\Ph4(q^2\pl 2q\pl 2)& q\Ph2\Ph4(q^2\pl 2q\pl 2)& (2q^2\pl 2q\pl 1)\Ph2^2\\

\phi_{20,11}& \frac12 q\Ph2\Ph6(2q^5\pl 3q^4\pl 5q^3\pl 4q^2\pl 4q\pl 2)& 
 q\Ph2\Ph4(q^2\pl 2q\pl 2)& q\Ph2\Ph4(q^2\pl 2q\pl 2)& (2q^2\pl 2q\pl 1)\Ph2^2\\

\phi_{20,10}& q(q^8\pl 2q^7\pl 2q^6\pl 3q^5\pl 3q^4\pl 4q^3\pl 2q^2\pl 2q\pl 1)& 
 q\Ph2\Ph4(q^2\pl 2q\pl 2)& q\Ph2\Ph4(q^2\pl 2q\pl 2)& (2q^2\pl 2q\pl 1)\Ph2^2\\

\phi_{15,12}'& q^2\Ph3\Ph6(q^3\pl q^2\pl 2q\pl 1)& 
 q\Ph3(2q^3\pl 2q\pl 1)& q\Ph3(2q^3\pl 2q\pl 1)& \Ph3(3q^2\pl q\pl 1)\\

\phi_{15,12}''& \frac12 q\Ph3\Ph6(2q^4\pl 3q^3\pl 3q^2\pl q\pl 1)& 
 q\Ph3(q^3\pl 2q^2\pl q\pl 1)& q\Ph3(q^3\pl 2q^2\pl q\pl 1)& q(2q\pl 3)\Ph3\\

\phi_{5,16}& \frac12 q(2q^8\pl q^7\pl q^6\pl 2q^4\pl 2q^3\pl q\pl 1)& 
 q^2(q^3\pl 2q^2\pl q\pl 1)& q^2(q^3\pl 2q^2\pl q\pl 1)&
 q(q^3\pl q^2\pl 2q\pl 1)\\

\phi_{10,14}& \frac12 q(3q^7\pl 3q^6\pl 4q^5\pl 4q^4\pl 2q^3\pl 4q^2\mn q\pl 1)& 
 q\Ph4(q^3\pl 2q^2\pl q\pl 1)& q\Ph4(q^3\pl 2q^2\pl q\pl 1)&
 q(q^3\pl 4q^2\pl 3q\pl 2)\\

\phi_{16,13}& q^2\Ph4(q^5\pl 2q^4\pl q^3\pl 2q^2\pl q\pl 1)& 
 q\Ph2\Ph4(2q^2\pl q\pl 1)& q\Ph2\Ph4(2q^2\pl q\pl 1)&  q\Ph2(3q^2\pl 2q\pl 3)\\

\phi_{16,15}& q^2\Ph4(q^5\pl 2q^4\pl q^3\pl 2q^2\pl q\pl 1)& 
q\Ph2\Ph4(2q^2\pl q\pl 1)& q\Ph2\Ph4(2q^2\pl q\pl 1)& 
q\Ph2(3q^2\pl 2q\pl 3)\\

\phi_{10,18}& q^3(2q^6\pl q^5\pl 2q^4\pl q^3\pl 2q^2\pl q\pl 1)& 
 q^2\Ph4(2q^2\pl 2q\pl 1)& q^2\Ph4(2q^2\pl 2q\pl 1)& 
 q(3q^3\pl 3q^2\pl 3q\pl 1)\\

\phi_{4,21}& \frac12 q^4\Ph2\Ph6(2q^2\pl q\pl 1)& q^3\Ph2\Ph4& q^3\Ph2\Ph4&
 q^2\Ph2^2\\

\phi_{4,23}&  \frac12 q^4\Ph2\Ph6(2q^2\pl q\pl 1)& q^3\Ph2\Ph4& q^3\Ph2\Ph4&
 q^2\Ph2^2\\

\phi_{4,24}& q^4(q^5\pl q^4\pl q^2\pl 1)& q^3\Ph2\Ph4& q^3\Ph2\Ph4& q^2\Ph2^2\\

\phi_{1,40}&  q^9&  q^6&  q^6&  q^4 \\   

\end{array}$

\newpage
$\begin{array}{c|ccccc}
 
 \# \cF\mn cl & \frac{1}{32}(l\mn 1)(l\mn 5)& \frac{1}{8}(l\mn 3)(l\mn 5) & \frac{1}{16}(l\mn 5)^2 & \frac{1}{96}(l\mn 3)(l\mn 5)(l\mn 13) & \frac{1}{7680}(l\mn 5)(l^3\mn 35l^2\pl 355l\mn 705) \\  
 C_\GG(s)(q) & B_2(q)\Phi_1^2 & A_2(q)\Phi_1^2 & A_1(q)^2\Phi_1^2 & A_1(q)\Phi_1^3 & \Phi_1^4\\
\hline
 
\phi_{1,0} & 1 & 1 & 1 & 1 & 1\\    

\phi_{4,1} & \frac12 q^3\pl q^2\pl \frac12 q\pl 2& q^2\pl q\pl 2& 2q\pl 2&
 q\pl 3& 4\\   

\phi_{4,3}& \frac12 q^3\pl q^2\pl \frac12 q\pl 2& q^2\pl q\pl 2& 2q\pl 2&
 q\pl 3& 4\\   

\phi_{4,4}& \Ph2(q^2\mn q\pl 2)& q^2\pl q\pl 2& 2q\pl 2& q\pl 3& 4\\   

\phi_{10,2}& 2q^3\pl 2q^2\pl 2q\pl 4& 3q^2\pl 3q\pl 4& q^2\pl 4q\pl 5&
 3q\pl 7& 10\\   

\phi_{16,3} & 4\Ph2\Ph4& q^3\pl 5q^2\pl 5q\pl 5& 2(q\pl 3)\Ph2&
 6q\pl 10& 16\\   

\phi_{16,5}& 4\Ph2\Ph4& q^3\pl 5q^2\pl 5q\pl 5& 2(q\pl 3)\Ph2& 6q\pl 10& 16\\

\phi_{15,4}''& q^4\pl 3q^3\pl 4q^2\pl 3q\pl 4& (q\pl 4)\Ph3&
 2q^2\pl 8q\pl 5& 6q\pl 9& 15\\

\phi_{5,8}& q^4\pl q^3\pl q\pl 2& 2q^2\pl 2q\pl 1& q^2\pl 2q\pl 2& 2q\pl 3& 5\\

\phi_{10,6}& (3q\pl 2)\Ph4& q^3\pl 3q^2\pl 3q\pl 3& q^2\pl 6q\pl 3&
 4q\pl 6& 10\\

\phi_{15,4}'& 4q^3\pl 4q^2\pl 4q\pl 3& (q\pl 4)\Ph3& 3q^2\pl 6q\pl 6&
 6q\pl 9& 15\\

\phi_{20,5}& 2q^4\pl \frac92 q^3\pl 5q^2\pl \frac92 q\pl 4& 
2q^3\pl 7q^2\pl 7q\pl 4& 2(2q\pl 3)\Ph2& 9q\pl 11& 20\\

\phi_{20,7}& 2q^4\pl \frac92 q^3\pl 5q^2\pl \frac92 q\pl 4& 
2q^3\pl 7q^2\pl 7q\pl 4& 2(2q\pl 3)\Ph2& 9q\pl 11& 20\\

\phi_{20,6}& \Ph2(2q^3\pl 3q^2\pl q\pl 4)& 2q^3\pl 7q^2\pl 7q\pl 4&
 2(2q\pl 3)\Ph2& 9q\pl 11& 20\\

\phi_{24,6}& 2\Ph3\Ph2^2& 4\Ph2\Ph3& 6\Ph2^2& 12q\pl 12& 24\\

\phi_{6,10}' & 2q\Ph3& \Ph2\Ph3& 2q^2\pl 2q\pl 2& 3\Ph2& 6\\

\phi_{6,10}''& 2q\Ph3& \Ph2\Ph3& 2q^2\pl 2q\pl 2& 3\Ph2& 6\\

\phi_{6,10}'''& \Ph3\Ph4& \Ph2\Ph3& q^2\pl 4q\pl 1& 3\Ph2& 6\\

\phi_{6,10}''''& \Ph3\Ph4& \Ph2\Ph3& q^2\pl 4q\pl 1& 3\Ph2& 6\\

\phi_{24,9}& 2q^4\pl 7q^3\pl 6q^2\pl 7q\pl 2& 4\Ph2\Ph3& 6\Ph2^2&
 12q\pl 12& 24\\

\phi_{24,7}& 2q^4\pl 7q^3\pl 6q^2\pl 7q\pl 2& 4\Ph2\Ph3& 6\Ph2^2&
 12q\pl 12& 24\\

\phi_{30,8}&  3q^4\pl 8q^3\pl 8q^2\pl 8q\pl 3& 5\Ph2\Ph3& 7q^2\pl 16q\pl 7&
 15q\pl 15& 30\\

\phi_{6,12}&  q^4\pl 2q^3\pl 2q\pl 1& \Ph2\Ph3& q^2\pl 4q\pl 1& 3\Ph2& 6\\

G_{29}[1]&  0& 0& 0& 0& 0\\

G_{29}[\ze_5]&  0& 0& 0& 0& 0\\

G_{29}[\ze_5^2]&  0& 0& 0& 0& 0\\

G_{29}[\ze_5^3]&  0& 0& 0& 0& 0\\

G_{29}[\ze_5^4]&  0& 0& 0& 0& 0\\

\phi_{20,9}&  
4q^4\pl \frac92 q^3\pl 5q^2\pl \frac92 q\pl 2& 4q^3\pl 7q^2\pl 7q\pl 2&
 2\Ph2(3q\pl 2)& 11q\pl 9& 20\\

\phi_{20,11}& 4q^4\pl \frac92 q^3\pl 5q^2\pl \frac92 q\pl 2& 
4q^3\pl 7q^2\pl 7q\pl 2& 2\Ph2(3q\pl 2)& 11q\pl 9& 20\\

\phi_{20,10}& \Ph2(4q^3\pl q^2\pl 3q\pl 2)& 4q^3\pl 7q^2\pl 7q\pl 2&
 2\Ph2(3q\pl 2)& 11q\pl 9& 20\\

\phi_{15,12}'& q(3q^3\pl 4q^2\pl 4q\pl 4)& (4q\pl 1)\Ph3& 6q^2\pl 6q\pl 3&
 9q\pl 6& 15\\

\phi_{15,12}''& 4q^4\pl 3q^3\pl 4q^2\pl 3q\pl 1& (4q\pl 1)\Ph3&
 5q^2\pl 8q\pl 2& 9q\pl 6& 15\\

\phi_{5,16}& 2q^4\pl q^3\pl q\pl 1& q(q^2\pl 2q\pl 2)& 2q^2\pl 2q\pl 1&
 3q\pl 2& 5\\

\phi_{10,14}& q(2q\pl 3)\Ph4& 3q^3\pl 3q^2\pl 3q\pl 1& 3q^2\pl 6q\pl 1&
 6q\pl 4& 10\\

\phi_{16,13}& 4q\Ph2\Ph4& 5q^3\pl 5q^2\pl 5q\pl 1& 2\Ph2(3q\pl 1)&
 10q\pl 6& 16\\

\phi_{16,15}& 4q\Ph2\Ph4& 5q^3\pl 5q^2\pl 5q\pl 1& 2\Ph2(3q\pl 1)&
 10q\pl 6& 16\\

\phi_{10,18}& 2q(2q^3\pl q^2\pl q\pl 1)& q(4q^2\pl 3q\pl 3)& 5q^2\pl 4q\pl 1&
 7q\pl 3& 10\\

\phi_{4,21}& \frac12 q(4q^3\pl q^2\pl 2q\pl 1)& q(2q^2\pl q\pl 1)& 2q\Ph2&
 3q\pl 1& 4\\

\phi_{4,23}&  \frac12 q(4q^3\pl q^2\pl 2q\pl 1)& q(2q^2\pl q\pl 1)& 2q\Ph2&
 3q\pl 1& 4\\

\phi_{4,24}&  q\Ph2(2q^2\mn q\pl 1)& q(2q^2\pl q\pl 1)& 2q\Ph2& 3q\pl 1& 4\\

\phi_{1,40}&  q^4&  q^3&  q^2&  q&  1\\   

\end{array}$

\normalsize


\end{document}